\documentclass[11pt]{article}
\usepackage{multicol}
\usepackage{amssymb,amsmath,amsthm,bm}
\usepackage[colorlinks=true, urlcolor=blue,linkcolor=blue, citecolor=blue]{hyperref}
\usepackage{mathrsfs,verbatim}
\usepackage{caption,cleveref }
\usepackage{graphicx, float}
\usepackage{color, mathtools, tikz, xcolor}
\usepackage{enumerate,enumitem}
\usepackage{tikz}
\usepackage[plain,lined,shortend,linesnumbered]{algorithm2e}
\usetikzlibrary{shapes.geometric}
\usepackage{pdfpages}
\usepackage{bbm}
\usetikzlibrary{calc}

\newcommand{\E}{\mathbb{E}}

\usepackage[mathscr]{euscript}
\usepackage{appendix}
\usepackage{cancel}
\usepackage[normalem]{ulem}
\usepackage{geometry}
\usepackage[affil-it]{authblk}
\numberwithin{equation}{section}
\newcommand{\cmt}[1]{} 
\newcommand{\F}{\mathbb{F}}

\renewcommand{\S}{\mathcal{S}}

\renewcommand{\l}{\left}
\renewcommand{\r}{\right}
\renewcommand{\sp}[1]{\mathrm{span}(#1)}

\usepackage{svg}
\usepackage{multicol}
\usepackage{amssymb,amsmath,amsthm,bm}
\usepackage[colorlinks=true, urlcolor=blue,linkcolor=blue, citecolor=blue]{hyperref}
\usepackage{pgfplots}

\numberwithin{equation}{section}

	\newtheorem{theorem}{Theorem}[section]

	\newtheorem{conjecture}[theorem]{Conjecture}
	\newtheorem{claim}{Claim}
	
	\newtheorem{lemma}[theorem]{Lemma}

        \theoremstyle{definition}
	
\newenvironment{proofclaim}[1][Proof of claim]{\begin{proof}[#1]}{\end{proof}}

\setlist{nolistsep}

\title{Maximum in-general-position set in a random subset of $\mathbb{F}^d_q$}

\author[2]{Yaobin Chen\thanks{Email: {\tt ybchen21@m.fudan.edu.cn}.
Supported by National Natural Science Foundation of China grant 123B2012.}}

\author[1]{Jiaxi Nie\thanks{Email: {\tt jnie47@gatech.edu}.}}

\author[2]{Jing Yu\thanks{Email: {\tt jyu@fudan.edu.cn}.} }

\author[2]{Wentao Zhang\thanks{Email: {\tt wtzhang20@fudan.edu.cn}.}}

\affil[1]{School of Mathematics, Georgia Institute of Technology, Atlanta, GA 30332, USA}

\affil[2]{Shanghai Center for Mathematical Sciences, Fudan University, Shanghai, 200438, China}

\begin{document}

\maketitle

\begin{abstract}
Let $\alpha(\mathbb{F}_q^{d},p)$ be the maximum possible size of a point set in general position in a $p$-random subset of $\mathbb{F}_q^d$. We determine the order of magnitude of $\alpha(\mathbb{F}_q^{d},p)$ up to a polylogarithmic factor by proving the balanced supersaturation conjecture of Balogh and Luo. Our result also resolves a conjecture implicitly posed by the first author, Liu, the second author and Zeng. In the course of our proof, we establish a lemma that demonstrates a ``structure vs. randomness'' phenomenon for point sets in finite-field linear spaces, which may be of independent interest.
\end{abstract}

\section{Introduction}
Let $d$ be a positive integer and $\mathbb{F}$ be a field. For $0\le k\le d$, a \textit{$k$-flat} of $\F^d$ is a $k$-dimensional affine subspace. A point set $P$ in $\mathbb{F}^{d}$ is in \textit{general position} if no $k+2$ points of $P$ lie in a $k$-flat for any $1\le k\le d-1$.
Sets in general position arise naturally in various central problems of discrete geometry and combinatorics.

A classical example is the famous \textit{no-three-in-line problem} raised by Dudeney \cite{Dudeney} in 1917, which asks if there exists a general position set in $\mathbb{R}^2$ containing $2n$ points from the grid $[n]\times [n]$. This question and its variations have been extensively studied for a long time, see \cite{hall1975some,flammenkamp1998progress,por2007no,lefmann2008no,SZ} for some latest results. 

 Another closely related problem, posed by Erd\H{o}s~\cite{erdos1986some} in 1986, is to determine the maximum size of a general position set contained in an arbitrary point set of size $n$ in $\mathbb{R}^2$ without collinear quadruples. Using the hypergraph container method, Balogh and Solymosi \cite{BS} showed that the answer is at most $n^{5/6+o(1)}$, which is a breakthrough result.

In this work, we focus on a random Tur{\'a}n-type problem related to general position sets in $\mathbb{F}^d_q$, where $\mathbb{F}_{q}$ is the finite field whose order is the prime power $q$. Random Tur\'an-type problems form an active area in extremal and probabilistic combinatorics concerning the maximum value of certain parameters in random structures. For example, see~\cite{kohayakawa1997k,conlon2016combinatorial,schacht2016extremal,jiang2022balanced,spiro2022random,Mubayi2023OnTR,nenadov2024number,nie2023tur,nie2023random}. 
More precisely, we study the maximum size of a general position set contained in a random subset of $\mathbb{F}_{q}^{d}$, where $d\ge 2$ is a fixed integer and $q$ is a sufficiently large prime power. This problem was initially considered by Roche-Newton and Warren~\cite{roche2022arcs} and Bhowmick~\cite{bhowmick2022counting}.
Given $p\in [0,1]$, let $\alpha(\mathbb{F}_q^d,p)$ denote the maximum size of a general position set contained in a \textit{$p$-random} set $\mathbf{S}_p \subset \mathbb{F}_q^d$, where each point is chosen independently with probability $p$.

Note that determining $\alpha(\mathbb{F}_q^d,1)$ reduces to the following fundamental question: 
\begin{center}
    How large can a general position set in $\mathbb{F}^d_q$ be? 
\end{center}

Erd\H{o}s~\cite{roth1951problem} observed that the \textit{moment curve} \[\{(x, x^2, \dots, x^d) \in \mathbb{F}_{q}^d : x \in \mathbb{F}_{q}\}\] is in general position, giving the lower bound $\alpha(\mathbb{F}_q^d,1) \ge q$. On the other hand, since $\mathbb{F}_q^d$ is a disjoint union of $q$ hyperplanes and each can contribute at most $d$ points to a general position set, $\alpha(\mathbb{F}_q^d,1) \le dq$.  Hence, we have $\alpha(\mathbb{F}_q^d,1) = O(q)$.

\begin{figure}[H]
    \centering
    \begin{tikzpicture}
        \newcommand {\xl} {5.5}
        \newcommand {\xo} {0.05}
        \newcommand {\xt} {1}
        \newcommand {\xth} {3}
        \newcommand {\xf} {5.2}
        \newcommand {\yo} {0.01732*5}
        \newcommand {\yl} {3.33}
        \newcommand {\ys} {1.75}
        \newcommand {\yt } {1.75}
        \newcommand {\yth} {1.75}
        \newcommand {\yf} {3}

        \begin{axis}[
        width = 0.7\textwidth,
        height = 0.4\textwidth,
        xlabel = {$p$},
        ylabel = {$\alpha(\mathbb{F}^d_q, p)$},
        xmin = 0, xmax = \xl,
        ymin = 0, ymax = \yl,
        xtick = {\xo, \xt, \xth, \xf},
        xticklabels = {$q^{-d}$, $q^{-d+1/d}$, $q^{-1+1/d+o(1)}$, $1$},
        ytick = {\ys, \yf},
        yticklabels = {$q^{1/d+o(1)}$, $q$},
        axis lines=middle,
        xlabel style={at=(current axis.right of origin), anchor=west},
        ylabel style={at=(current axis.above origin), anchor=south},
        ]
        \addplot[color = black]
            plot coordinates {
                (\xo,0)
                (\xt,\yt)
                (\xth,\yth)
                (\xf,\yf)
            };
    
        \addplot[color = gray,dashed]
            plot coordinates {
                (\xt, 0)
                (\xt, \yt)
            };
        \addplot[color = gray,dashed]
            plot coordinates {
                (0, \ys)
                (\xt, \ys)
            };
        \addplot[color = gray,dashed]
            plot coordinates {
                (\xth, 0)
                (\xth, \yth)
            };
        \addplot[color = gray,dashed]
            plot coordinates {
                (\xf, 0)
                (\xf, \yf)
            };
        \addplot[color = gray,dashed]
            plot coordinates {
                (0, \yf)
                (\xf, \yf)
            };
        \end{axis}

    \end{tikzpicture}
    \caption{Behavior of $\alpha(\mathbb{F}_q^d,p)$.}
    \label{fig:randomturan}
    \end{figure}
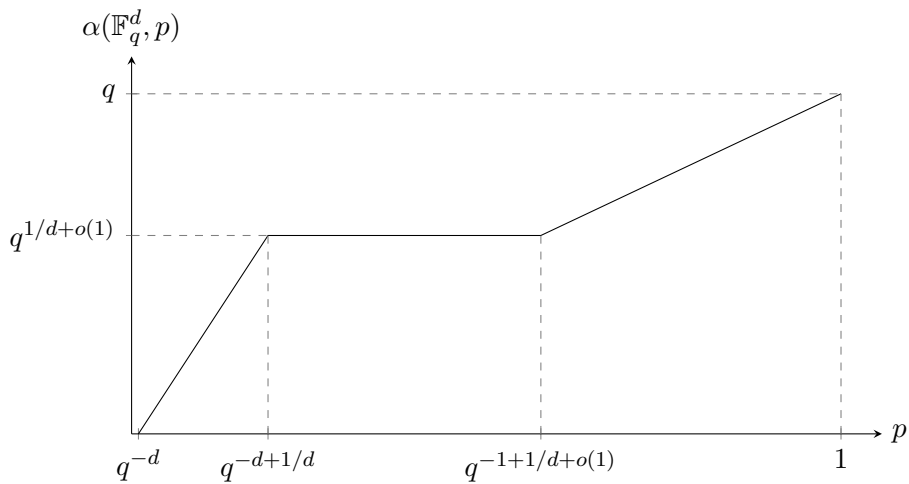

For $d=2$, Roche-Newton--Warren~\cite{roche2022arcs} and Bhowmick--Roche-Newton~\cite{bhowmick2022counting} established essentially tight bounds for $\alpha(\mathbb{F}_q^2,p)$ when $p\le q^{-1+o(1)}$ or $p\ge q^{-1/3+o(1)}$. In~\cite{CLNZ}, the first author, Liu, the second author and Zeng determines the order of magnitude of $\alpha(\mathbb{F}_q^2,p)$ up to polylogarithmic factors for all possible values of $p$ and implicitly posed  a conjecture about the order of magnitude of $\alpha(\mathbb{F}_q^d,p)$ for all fixed integer $d\ge 3$ as in Figure~\ref{fig:randomturan}. Recently, Balogh and Luo~\cite{LB} confirms the case $d=3$. These previous works are all grounded in the hypergraph container method developed independently by Balogh--Morris--Samotij~\cite{balogh2015independent} and Saxton--Thomason~\cite{saxton2015hypergraph} together with suitable balanced supersaturation results.

\textit{Balanced supersaturation} refers to the phenomenon that a sufficiently ``dense'' structure must contain a collection of evenly distributed objective substructures. Establishing suitable balanced supersaturation is usually the main challenge in random Turán-type problems, see for example, \cite{nie2023random,nie2023tur, Mubayi2023OnTR,MorrisC2l, jiang2022balanced,nie2024random2, mckinley2023random}.

A \textit{coplanar $k$-set} is a set of $k$ points in a $(k-2)$-flat, i.e., its affine span has dimension $k-2$.  Thus a point set in general position is a point set with no coplanar sets. In~\cite{LB}, Balogh and Luo posed a balanced supersaturation conjecture for the coplanar sets in $\mathbb{F}_q^d$ which, together with the hypergraph container method, would imply the order of magnitude of $\alpha(\mathbb{F}_q^d,p)$. 

In this paper, we confirm their balanced supersaturation conjecture.

\begin{theorem}\label{main}
     For any integer $d\ge 2$, there exists a constant $C>0$ such that for any sufficiently large prime power $q$ the following holds. If $X\subseteq \F^d_q$ with $|X|=n\ge Cq$, then there is a collection $\mathcal{S}$ of coplanar $(d+1)$-sets in $X$ such that 
    \begin{enumerate}[label=\rm{(\roman*)}]
        \item\label{equationtotal1} $|\mathcal{S}|\ge \frac{1}{C} n^{d+1}q^{-1}$.
        \item\label{equationtotal2} $\Delta_j(\mathcal{S})\le Cn^{d+1-j}\cdot{q^{-\frac{d+1-j}{d}}}$ for $1\le j\le d+1$. 
    \end{enumerate}
\end{theorem}

Hence, by Claim~3 in \cite{LB}, we obtain an affirmative resolution of the conjecture of the first author, Liu, the second author and Zeng. An illustration of this phenomenon is provided in Figure~\ref{fig:randomturan}.

\begin{theorem}\label{thm:randomturan1}
As the prime power $q\to \infty$, asymptotically almost surely, we have 
\begin{equation*}
\alpha(\mathbb{F}_q^d,p)=
\l\{
    \begin{aligned}
    &\Theta(pq^{d}),~~~&q^{-d+o(1)}\le p\le q^{-d+1/d-o(1)},\\
    &q^{1/d+o(1)},~~~&q^{-d+1/d-o(1)}\le p\le q^{-1+1/d+o(1)},\\
    &\Theta(pq),~~~&q^{-1+1/d+o(1)}\le p\le 1,
\end{aligned}
\r.
\end{equation*}
where all $q^{o(1)}$ factors are polylogarithmic. 
\end{theorem}

 Our argument builds on the ``local probabilistic'' strategy employed by Balogh and Luo~\cite{LB}, incorporating a few novel ideas to address difficulties that arise only in higher dimensions. A central component is a key lemma that establishes a ``structure vs. randomness'' phenomenon. Its proof relies on extremal results on tight Hamiltonian cycles in uniform hypergraphs, classical hypergraph Tur{\'a}n theorems, and a related supersaturation lemma. This key lemma naturally splits our proof into the relatively simpler ``structure'' case and the more intricate ``randomness'' case. To handle the latter, we introduce an operation called support-enlarging substitution to ensure the critical coplanarity of a point set. A sketch of the argument is given in Section~\ref{proofoutline}, with full details provided in subsequent sections.

 \textbf{Organization of the paper.}
Section~\ref{pre} introduces some notations and collects some useful lemmas, which the reader may skip on a first reading and return to as needed. Section~\ref{proofoutline} provides a sketch of our overall proof. In Section~\ref{blanceprop}, we establish the key ``structure vs. randomness'' lemma, which divides the proof into the ``structure'' case and the ``randomness'' case. Section~\ref{sec:struc} resolves the ``structure'' case while Section~\ref{BH} and~\ref{BL} prove the ``randomness'' case. In Section~\ref{proof1} we combine these ingredients to prove Theorem~\ref{main}. Finally, Section~\ref{concluding} contains some concluding remarks.

\section{Preliminary}\label{pre}
\subsection{Notation}
\par In this section, we introduce notation and basic properties of flats in $\mathbb{F}_q^d$. 

For a flat $F\subseteq \mathbb{F}_q^d$, we denote its dimension by $\dim (F)$. A point set $X$ is \textit{IGP} if it is in general position; otherwise, it is \textit{coplanar}. We say $X$ is \textit{critical coplanar} if $X$ is coplanar but all of its proper subsets are IGP. Given $X\subseteq \mathbb{F}_q^d$, we write $F_X$ or $\sp{X}$ for the flat determined by $X$, i.e., the minimal flat containing $X$. Equivalently,
$$
F_X\;=\;\sp{X}\;=\;\l\{\sum_{x\in X}\lambda_x{x}:\ \sum_{x\in X}\lambda_x=1,~\lambda_x\in\F_q\r\}.
$$

We say that $U$ is \textit{supported} by $X$ if $U\subseteq F_X$ but $U\nsubseteq F_{X'}$ for any $X'\subsetneq X$; in this case we also call $X$ is a \textit{support} of $U$. 

Let $F\subseteq \mathbb{F}_q^d$ be a flat. For $s\le \dim F$ and $\varepsilon<1$, we say $F$ is \textit{$(s,\varepsilon, X)$-heavy}  if there exists an $s$-flat $F'\subseteq F$ such that 
$$|X\,\cap\, F'|\;\ge\; \varepsilon^{\dim (F)-s}\,|X\,\cap\, F|.$$

In particular, $F$ is \textit{$(s,\varepsilon, X)$-balanced} if it is $(s,\varepsilon, X)$-heavy but not $(s-1, \varepsilon, X)$-heavy. Furthermore, we say $F$ is \textit{$(s,j,\varepsilon, X)$-balanced} if it is $(s,\varepsilon, X)$-balanced and 
\[|X\,\cap\, F|\;\in\; (\varepsilon^{d-s}nq^{-\frac{j+1}{d}}, \varepsilon^{d-s}nq^{-\frac{j}{d}}].\] 

For $\l(d-1\r)\ge j>i\ge 1$, an IGP $(j+1)$-set $J$ is \textit{$(i,\varepsilon,X)$-good} if $|X\cap F_J|\le nq^{-\frac{j}{d}}$ and there exists a partition $V_1,\dots,V_{i+1}$ of $J$ with 
\[
|V_1| \;=\; \left\lceil \tfrac{j-i}{2} \right\rceil + 1, 
\qquad 
|V_2| \;=\; \left\lfloor \tfrac{j-i}{2} \right\rfloor + 1 
\quad \text{and }\quad
|V_3| \;=\; \cdots \;=\; |V_{i+1}| \;=\; 1 
\]
such that for any $(i+1)$-set $I$ with $|I\cap V_s|=1$
for $s\in [i+1]$, $F_I$ is $(i,j,\varepsilon ,X)$-balanced. For each such $(i,\varepsilon,X)$-good IGP $(j+1)$-set $J$, we arbitrarily designate a vertex in $V_1$ as the \textit{root} of $J$. {Let $\mathcal{J}_{i,j}$ be the collection of $(i,\varepsilon,X)$-good $(j+1)$-sets, and define the \textit{root function} $f:\mathcal{J}_{i,j}\to X$ that maps each $(i,\varepsilon,X)$-good $(j+1)$-set to its root.}

We say a bipartite graph $H$ with parts $A$ and $B$ is \textit{bi}-$(a,b)$ if $\deg_{B}(v)\ge a$ for every $v\in A$ and $\deg_{A}(v)\le b$ for every $v\in B$. For a graph $H$, a \textit{directed Hamiltonian path} $P$ of $H$ is an enumeration $v_1,\dots ,v_n$ of vertices in $V(H)$ such that $\overrightarrow{v_iv_{i+1}}\in E(H)$ for $i\in [n-1]$. We call $v_1$ the \textit{start point of $P$}. We call a hypergraph \textit{$k$-graph} if every edge contains exactly $k$ vertices. A \textit{tight Hamiltonian cycle} of a $k$-graph $G$ is a cyclic ordering $v_1,\dots, v_n$ of the vertices of $G$ such that $v_iv_{i+1} \dots v_{i+k-1}$ form an edge for every $1\le i\le n$, where indices are taken modulo $n$. For a $k$-graph $G$, the \textit{degree} of $v\in V(G)$, denoted by $\deg_{G}(v)$, is the number of edges containing $v$. We use $\delta(G)$ to denote the \textit{minimum vertex degree} of $G$, that is, $\delta(G)=\min_{v\in V(G)}\deg_G(v)$.

{In this paper, we use the notation $f(x)=O(g(x))~(f(x)=\Omega(g(x)))$ if there exists a constant $C>0$ such that $f(x)<Cg(x)~(f(x)>Cg(x),\text{respectively})$ as $x\rightarrow \infty$. 
Similarly, we write $f(x)=o(g(x))~(f(x)=\omega(g(x)))$ if for any constant $C>0$, $f(x)<Cg(x)~(f(x)>Cg(x),\text{respectively})$ as $x\rightarrow \infty$.}

{Throughout, we will omit floors and ceilings when they do not affect the argument. The
constants in the hierarchies used to state our results are chosen from right to left. For
example, if we state that a result holds whenever 
\[ 0 \;<\; 1/q \;\ll\; a \;\ll\; b \;\ll\; c \;\le\; 1 \]
(where $q$ is
the order of finite field $\mathbb{F}_{q}$), then there exists a non-decreasing function $f : (0, 1] \rightarrow (0, 1]$
such that the result holds for all $0 < a, b, c \le  1$ and all $q \in  \mathbb{N}$ with 
\[ b \;\le\; f(c), \; a \;\le\; f(b), \; 1/q \;\le\; f(a). \]
Hierarchies with more constants are defined in the same manner.}

\subsection{Basic properties and tools}
\par In our proofs, we will frequently use the following auxiliary results.

We begin with a supersaturation lemma for $k$-partite $k$-graphs, which will be used to construct coplanar point sets.

\begin{lemma}[Supersaturation Lemma {\cite{erdos1983supersaturated}}]\label{supersaturation}
    Let $F$ be a $k$-partite $k$-graph with $k\ge 2$. For any $\varepsilon>0$, there exist positive $\delta=\delta(\varepsilon, F)$ and $n_0=n_0(\varepsilon, F)$ such that for all $n\ge n_0$, every $n$-vertex $k$-graph with more than $\varepsilon n^k$ edges contains at least $\delta n^{v(F)}$ copies of $F$, where $v(F)=|V(F)|$.
\end{lemma}

Besides, we will frequently use the following version of Chernoff bound (see Corollary A.1.14 in \cite{alonprob}) to compute the probability of our probabilistic construction.

\begin{lemma}[Chernoff bound]\label{lem: Chernoff}
   Let $X_1,\dots, X_n$ be independent Bernoulli random variables, and let $S_n=\sum_{i=1}^n X_i$ with $\mu=\E(S_n)$. Then for all $\varepsilon>0$, 
\[
    \Pr\!\left[\,|S_n - \mu| \ge \varepsilon \mu\, \right] 
    \;\le\; \exp\!\left(-2c_\varepsilon \mu\right).
\]
where $c_{\varepsilon}=\min \{-\ln (e^{\varepsilon}(1+\varepsilon)^{-(1+\varepsilon)}), \varepsilon^2/2\}$. In particular, for $h\ge20\mu$, 
    \[
    \Pr[S_n>h]\;\le\; 2\exp(-h/2).
    \]
\end{lemma}
We also need two lemmas concerning the existence of Hamiltonian cycles and paths.
\begin{lemma}
[\cite{Langham}]\label{Lemma:Hamilton}
    Let $k\ge 2$ and let $G$ be a $k$-graph. If 
    \[\delta(G) \;\ge\; \l(1-\frac{1}{k-1}\r)\binom{n-1}{k-1},\] then $G$ contains a tight Hamiltonian cycle.
\end{lemma}

\begin{lemma}\label{prop:deirectedhamiltonianpath}
    Let $G$ be a complete bipartite graph with bipartition $(A,B)$, where $|A|+1\ge |B|\ge |A|$. Then for each vertex $v\in B$, there exists a directed Hamiltonian path starting at $v$ that alternates between $A$ and $B$.
\end{lemma}

\begin{proof}
    If $|A|=n$ and $|B|=n+1$, label $A$ as $a_1,\dots,a_n$, set $b_1 = v$ and label $B \setminus \{v\}$ as $b_2,\dots,b_{n+1}$. Since $G$ is complete bipartite, the directed path $b_1, a_1, b_2, a_2, \dots, a_n, b_{n+1}$ is the Hamiltonian path as required.
    
    If $|B|=n$, the construction is the same but we omit  $b_{n+1}$.
\end{proof}

The following basic properties of flats will be used frequently throughout the paper.
\begin{lemma}\label{propflat}
    Let $X\subseteq\mathbb{F}_q^d$ and $F_1,F_2$ be flats in $\mathbb{F}_q^d$. Then the following properties hold:
    \begin{enumerate}[label=\rm{(F\arabic*)}]

        \item\label{propflat2} If $X$ is IGP, then for any $w\in X$ and $w'\in F_X\setminus F_{X\setminus \{w\}}$, $\l(X\setminus \{w\}\r)\cup \{w'\}$ is IGP and \[F_{\l(X\setminus \{w\}\r)\,\cup\, \{w'\}}\;=\;F_{X}.\] 
        \item\label{propflat3} $F_1\cap F_2$ is a flat. Moreover, if $F_1\nsubseteq F_2$ and $F_2\nsubseteq F_1$, then \[\dim(F_1\cap F_2)\;<\;\min\l\{\dim(F_1),\dim(F_2)\r\}.\]
        \item\label{propflat4} If $X$ is IGP, then any nonempty subset of $X$ is also IGP.
        \item\label{propflat5} If $X$ is IGP and $Y$, $Z$ are two nonempty subsets of $X$ with $Y\not\subseteq Z$ and $Z\not\subseteq Y$, then $F_{Y}\not\subseteq F_{Z}$ and $F_{Z}\not\subseteq F_{Y}$.
    \end{enumerate}
\end{lemma}

\begin{proof}
    
    For \ref{propflat2}, let $X\setminus\{w\}=\{x_0,x_1,\dots,x_k\}$. Since $X$ is IGP and $w'\in F_X\setminus F_{X\setminus \{w\}}$, $w'-x_0$ has a unique expression $\lambda(w-x_0)+\lambda_1(x_1-x_0)+\dots+\lambda_k(x_k-x_0)$ where $\lambda\neq0$. So
    \[
w - x_0 
   \;=\; \dfrac{1}{\lambda} \,\Bigl[
      (w' - x_0) \;-\; \lambda_1(x_1 - x_0) \;-\; \cdots \;-\; \lambda_k(x_k - x_0)
   \Bigr].
\]
    For every $y\in F_X$, $y-x_0$ can be expressed as a linear combination of $w-x_0$, $x_1-x_0$, $\dots$, $x_k-x_0$, and hence also as a linear combination of $w'-x_0$, $x_1-x_0$, $\dots$, $x_k-x_0$, establishing the claim.

    For \ref{propflat3}, if $F_1\cap F_2\neq\emptyset$, then take $x\in F_1\cap F_2$ and write $F_1 = x+U_1, F_2 = x+U_2$ with $U_1$, $U_2$ subspaces of $\mathbb{F}_q^d$. Then $F_1\cap F_2 = x + (U_1\cap U_2)$, which is also a flat. Take $y\in F_1 \setminus F_2$, then $y-x\in U_1 \setminus U_2$, so $\dim(U_1)>\dim(U_1\cap U_2)$, implying $\dim(F_1\cap F_2)<\dim(F_1)$. The same holds with $F_1, F_2$ swapped.

    To see \ref{propflat4}, without loss of generality let $X=\{x_0, x_1,\dots,x_k\}$ and $X'=\{x_0, x_1,\dots,x_j\}$ with $j<k$. Since $X$ is IGP, the vector set $\{x_1-x_0,\dots,x_k-x_0\}$ is linear independent. As its subset, $\{x_1-x_0,\dots,x_j-x_0\}$ is also linear independent, which suggests $X'$ is IGP.

    To see \ref{propflat5}, suppose $F_{Y}\subseteq F_{Z}$. The case when either of $Y$ and $Z$ is single point is trivial. If $|Y|$, $|Z| \ge 2$, let $Y=\{y_0,y_1,\dots,y_k\}$ and $Z=\{z_0,z_1,\dots,z_j\}$. Choose $y_i \in Y\setminus Z$. Then $y_i-z_0$ lies in $\sp{\{z_1-z_0, \dots, z_j-z_0\}}$, so $\{y_i, z_0, z_1, \dots, z_j\}$ is coplanar, contradicting \ref{propflat4}.
\end{proof}

Below are some properties concerning ``supported'' and critical coplanar set. 
\begin{lemma}\label{exactprop}
    Let $U,Y$ be point sets of $\mathbb{F}_{q}^d$. Suppose that $U$ is supported by $X\subseteq Y$. Then the following properties hold:
    \begin{enumerate}[label=\rm{(E\arabic*)}]\label{exactly}
        \item\label{exactprop1} $U\subseteq F_{Y'}$ for every $Y'\subseteq Y$ containing $X$. Moreover, if $Y$ is IGP, then $U\nsubseteq F_{Y''}$ for any $Y''\subseteq Y$ not containing $X$. 
        \item\label{exactprop2} If $Y$ is IGP, then $X$ is unique in $Y$.
        \item\label{exactprop3} If $Y$ is IGP and $U = \{u\}$, then $\{u\}\cup X$ is a critical coplanar set. 
    \end{enumerate}
\end{lemma}

\begin{proof}
    For \ref{exactprop1}, since $Y'$ contains $X$, we have $U\subseteq F_X\subseteq F_{Y'}$. When $Y$ is IGP, for $Y''\subseteq Y$ with $X\nsubseteq Y''$, assume $U\subseteq F_{Y''}$. Since $U$ is supported by $X$, we also have $Y''\nsubseteq X$. By the fact that $Y$ is IGP and \ref{propflat5}, we have $F_X\nsubseteq F_{Y''}$ and $F_{Y''}\nsubseteq F_X$. By Lemma \ref{propflat} \ref{propflat3}, $F_X\cap F_{Y''}$ has dimension strictly less than $\min\{\dim(F_X),\dim(F_{Y''})\}$. Choose $u\in U$ outside $F_X\cap F_{Y''}$, then $u$ lies in  both $F_{Y''}\setminus F_X$ and $F_X\setminus F_{Y''}$, contradicting the assumption that $X$ covers $U$.
    
    \ref{exactprop2} is a direct corollary of \ref{exactprop1}.
    
    To see \ref{exactprop3}, let $X=\{x_0,x_1,\dots,x_k\}$, it suffices to show $(X\setminus\{x_i\})\,\cup\,\{u\}$ is IGP for any $0\le i\le k$. Suppose not,  there exists some $i$ such that $F_{(X\setminus\{x_i\})\,\cup\,\{u\}}$ has dimension at most $|X|-2$. By Lemma \ref{propflat} \ref{propflat4}, $X\setminus\{x_i\}$ is an IGP set of $|X|-1$ points, so $F_{(X\setminus\{x_i\})\,\cup\,\{u\}}$ has dimension exactly $|X|-2$. This forces $F_{(X\setminus\{x_i\})\,\cup\,\{u\}} = F_{X\setminus\{x_i\}}$, so $u\in F_{X\setminus\{x_i\}}$, contradicting the fact that $u$ is supported by~$X$.
\end{proof}

Next, we present a lemma about enlarging the size of a support.
\begin{lemma}
    Let $X$ be an IGP set in $\mathbb{F}_q^d$ and $u$ be a point supported by $U\subsetneq X$. For any $w\in U$ and any subset $Y\subseteq X$ with $w\in Y$, the following hold:
    \begin{enumerate}[label=\rm{(C\arabic*)}]\label{enlargecoverset}
        \item\label{enlargecoverset1} $F_{Y}\cap F_{X\setminus\{w\}}=F_{Y \setminus\{w\}}$; in particular, $\dim(F_{Y}\cap F_{X \setminus \{w\}})=\dim(F_Y)-1=|Y|-2$.
        \item\label{enlargecoverset2} For each $y\in U\cup Y\setminus\{w\}$, $\dim(F_{U \cup Y \,\cup\, \{u\}\setminus \{w,y\}}\cap F_Y)<\dim(F_Y)$. 
        \item\label{enlargecoverset3} If there exists a point $$w'\in\; F_{Y}\setminus \left(\bigcup_{y\in U\cup Y}F_{U\cup Y \,\cup\, \{u\}\setminus\{w,y\}}\cup F_{X\setminus\{w\}}\right),$$ then $u$ is supported by $(U\cup Y\cup \{w'\})\setminus\{w\}$.  
    \end{enumerate}
\end{lemma}
\begin{proof}

        For~\ref{enlargecoverset1}, note that $Y \setminus\{w\}=Y \cap (X\setminus \{w\})$. This implies $F_{Y\setminus\{w\}}\subseteq F_Y \cap F_{X\setminus\{w\}}$. By Lemma~\ref{propflat}~\ref{propflat3}, $F_Y \,\cap\, F_{X\setminus\{w\}}$ is a flat in $F_Y$. If $w \in F_{X \setminus\{w\}}$, then $w\in \mathrm{span}(X \setminus\{w\})$, contradicting that $X$ is IGP. Thus $F_Y \cap F_{X\setminus\{w\}}\subseteq F_{Y\setminus\{w\}}$, proving equality.

For~\ref{enlargecoverset2}, it suffices to exhibit a point of $F_Y$ not contained in $F_{U\cup Y\,\cup\,\{u\}\setminus\{w,y\}}$. If $y\in U$, assume for a contradiction that $w\in F_{U\cup Y\,\cup\,\{u\}\setminus\{w,y\}}$. Since $X$ is IGP, $w\notin \mathrm{span}(U\cup Y\setminus\{w,y\})$. Hence there exist coefficients with sum $1$ such that
    \[
      w\;=\sum_{z\in U\cup Y\setminus\{w,y\}}\lambda_z z+\lambda_u u,
      \qquad \lambda_u\ne 0.
    \]
    As $u$ is supported by $U$, write $u=\sum_{z\in U}\xi_z z$ with $\sum_{z\in U}\xi_z=1$ and all $\xi_z\ne0$, which implies $w\in \mathrm{span}(U\cup Y\setminus\{w\})\subseteq \mathrm{span}(X\setminus\{w\})$, a contradiction. If $y\in Y\setminus U$ and $y\in F_{U\cup Y\,\cup\,\{u\}\setminus\{w,y\}}$, then using $u\in \mathrm{span}(U)$, we would get $\text{span}(Y\cup U\cup\{u\}\setminus \{y\})=\text{span}(Y\cup U\setminus \{y\})$. Hence, $y\in \text{span}(Y\cup U\cup\{u\}\setminus \{w,y\})\subseteq \text{span}(Y\cup U\setminus \{y\})\subseteq \text{span}(X\setminus y)$, contradicting IGP.

   Thus~\ref{enlargecoverset2}~holds.

 For~\ref{enlargecoverset3}, we must show $u\in F_{(U\cup Y\cup\{w'\})\setminus\{w\}}$ and and the minimality of $(U\cup Y\cup\{w'\})\setminus\{w\}$. Since $w'\in F_Y\setminus F_{X\setminus\{w\}}$, by~\ref{enlargecoverset1} we have $w'\notin F_{Y\setminus\{w\}}$, hence $F_{(Y\setminus\{w\})\cup\{w'\}}=F_Y$. Thus
    \[
      u\;\in\; \mathrm{span}(U)\;\subseteq\; \mathrm{span}\big((U\setminus\{w\})\cup (Y\setminus\{w\})\cup\{w'\}\big).\]
      If $u\in F_{Y'}$ for some strict subset $Y'\subsetneq (U\cup Y\cup\{w'\})\setminus\{w\}$, 
    let $x$ be a point in $(U\cup Y\cup \{w'\})\setminus \{w\}$ not in $Y'$. Then there is a linear combination of $u$ whose coefficients sum to 1 such that 
    \begin{align*}
        {u}\;=\sum_{z\in (U\cup Y\cup \{w'\})\setminus \{w,x\}}\lambda_{z} {z}\;=\sum_{z\in (U\cup Y)\setminus \{w,x\}}\lambda_{z} {z}\;+\;\lambda_{w'}{w'}.
    \end{align*}
    Note that $\lambda_{w'}\neq 0$, because otherwise $u\in \text{span}((U\cup Y)\setminus \{w\})$. By~Lemma~\ref{exactprop}~\ref{exactprop1}, since $U\nsubseteq (U\cup Y)\setminus \{w\}$, we have $u\notin \text{span}((U\cup Y)\setminus \{w\})$, a contradiction. So $\lambda_{w'}\neq 0$ indeed. Hence, $w'\in \text{span}(\l(U\cup Y\cup \{u\}\r)\setminus\{w,x\})$, contradicting the assumption that $w'\notin F_{\l(U\cup Y\,\cup\, \{u\}\r)\setminus\{w,x\}}$.
    Thus~\ref{enlargecoverset3} holds. 
\end{proof}

Next, we show a useful lemma about $(i,\varepsilon,X)$-good sets.

 \begin{lemma}\label{key ob}
     Let $i\in [d]$ and $\varepsilon\le 1/2$. Let $X$ be a point set of $\mathbb{F}_q^d$ and $J$ be an $(i,\varepsilon,X)$-good $(j+1)$-set $J=\{y,x_2,\dots,x_{j+1}\}$ with partition $V_1,\dots, V_{i+1}$ and $y\in V_1$. Then for any $(j-1)$-flat $F\subseteq F_J\setminus \{y\}$, we have 
\[
   \frac{|F_J \cap X|}
        {\;\bigl|(F_J \setminus F)\cap X\bigr|}
   \;\le\;
   \varepsilon^{-d}\, q^{1/d}.
\]
 \end{lemma}
\begin{proof}
    Since $J$ is an $(i,\varepsilon,X)$-good $(j+1)$-set, we have $|F_J\cap X|\le nq^{-\frac{j}{d}}$ and there exists an $(i+1)$-set $I$ containing $y$ with $F_I$ $(i,j,\varepsilon,X)$-balanced. Note that $y\notin F$. This implies that $F_I\cap F$ is a flat in $F_I$ with dimension at most $(i-1)$, and hence $|\l(F_{I}\setminus F\r)\cap X|\ge (1-\varepsilon)\varepsilon^{d-i}nq^{-\frac{j+1}{d}}$. Since $F_{I}\setminus F\subseteq F_{J}\setminus F$, we can see that 
\[
   \frac{|F_J \cap X|}
        {\;\bigl|(F_J \setminus F)\cap X\bigr|}\;\le\; \frac{|F_J \cap X|}
        {\;\bigl|(F_I \setminus F)\cap X\bigr|}\;\le\; \frac{n\,q^{-j/d}}
        {(1-\varepsilon)\,\varepsilon^{\,d-i}\,n\,q^{-(j+1)/d}}\;\le\;\varepsilon^{-d}\,q^{1/d}.\qedhere\]
\end{proof}

A key step in the proof of Theorem~\ref{main} is to map a coplanar set to a special one with desirable properties. The construction of this mapping relies on the following lemma.

\begin{lemma}[Auxiliary Graph Lemma]\label{auxiliary graph}
         For any integer $d\ge 2$, there exists constants $C=C(d)$ and $\varepsilon=\varepsilon(d)\in (0, 1)$ such that the following holds.
        
Let $F$ be a flat in $\mathbb{F}_q^d$. Suppose $X \subseteq \mathbb{F}_q^d$ satisfies $|X\cap F|\ge C$, and $F$ is $(s,\varepsilon, X)$-balanced with some $s\le \dim(F)$, witnessed by an $s$-flat $F'$. Then there exists a bi-$(3d^2, 6d^2\tfrac{|X\cap F|}{|X\cap F'|})$ graph $H_F$ with parts $A=X \cap F, B=X\cap F'$ such that for every $v \in A$, its neighborhood $N_{H_F}(v)$ cannot be covered by $2d$ distinct~$(s-1)$-flats.         
\end{lemma}
\begin{proof}
    We first construct an auxiliary $(2sd+1)$-graph $G$ with vertex set $X\cap F'$. A $(2sd+1)$-subset $\{x_1, \dots, x_{2sd+1}\} \subseteq X \cap F'$ is an edge of $G$ if it can be ordered so that:
\begin{itemize}
 \item $\{x_1, \dots, x_s\}$ is IGP, and
    \item for every $s \le i \le 2sd+1$, $x_i$ does not lie in the flat determined by any $s$ earlier vertices.
\end{itemize}
In particular, no $s+1$ vertices in such an ordering lie in the same $(s-1)$-flat.

By the $(s,\varepsilon,X)$-balanced property, we have 
\[
   |X \cap F'| \;\ge\; \varepsilon^{\,i-s}\,|X \cap F|
   \;\ge\; \varepsilon^{\,i-s} C,
\]
and every $(s-1)$-flat $F'' \subseteq F$ satisfies
\[
   |X \cap F''|
   \;\le\; \varepsilon^{\,i-s+1}\,|X \cap F|
   \;\le\; \varepsilon\,|X \cap F'|.
\]

Fix $x_1 \in X\cap F'$. We iteratively choose $x_2,\dots,x_{2sd+1}$, each time avoiding any $(s-1)$-flat spanned by $s$ earlier vertices. Each such forbidden flat contains at most $\varepsilon |X\cap F'|$ points, and there are at most $2^{2sd}$ such flats at each step. Thus each choice has at least $(1 - 2^{2sd} \varepsilon)|X\cap F'|$ options. Consequently, $x_1$ lies in at least
\[
   \frac{(1 - 2^{2sd}\,\varepsilon)^{2sd}\,|X \cap F'|^{\,2sd}}
        {(2sd)!}
\]
edges of $G$.

If $C$ and $1/\varepsilon$ are sufficiently large (depending on $d$), Lemma~\ref{Lemma:Hamilton} guarantees a tight Hamiltonian cycle $HC$ in $G$. Let $HC$ have edges $e_1, \dots, e_m$. Partition $A = X\cap F$ equitably into $A_1, \dots, A_m$ with $\bigl||A_i| - |A_j|\bigr| \le 1$.

We now define $H_F$: connect $u \in A_i$ to all $v \in e_i$ for each $i$. Then:
\begin{itemize}
    \item Every $u \in A$ has degree $|e_i| = 2sd+1 \le 3d^2$.
    \item If $N_{H_F}(u)$ could be covered by $2d$ distinct $(s-1)$-flats, then by the pigeonhole principle, some flat would contain at least $s+1$ neighbors of $u$, contradicting the construction of $G$.
    \item Each $v \in B$ lies in at most $3d^2$ edges of $HC$, so its degree in $H_F$ is at most
    \[
   3d^2 \,\Biggl\lceil \frac{|X \cap F|}{m} \Biggr\rceil
   \;\le\;
   6d^2 \cdot \frac{|X \cap F|}{|X \cap F'|}.
\]
\end{itemize}

Thus $H_F$ is bi-$(3d^2,\, 6d^2 \tfrac{|X\cap F|}{|X\cap F'|})$ with the claimed covering property.
\end{proof}

\section{Proof outline}\label{proofoutline}

\par The proof of Theorem~\ref{main} can be divided into two parts. 

In the first part, we prove a lemma demonstrating a ``structure vs. randomness'' phenomenon: given a point set $X\subseteq\F^d_q$ with $|X|=n\gg q$, either it has some nice structures -- there are $\Omega(n^{d+1})$ coplanar $(d+1)$-sets; or it has some pseudorandom properties -- for some $1\le i\le d-1$, there are $\Omega(n^{i+1})$ IGP $(i+1)$-sets $I$ such that $F_I$ contains reasonably many points in $X$ and they are evenly distributed, i.e., $F_I$ is ``dense'' and ``balanced''. In the ``structure'' case, it suffices to sample each coplanar $(d+1)$-set independently with a suitable probability to obtain a family with the required size and degree bounds. The ``randomness'' case, however, requires much more technical discussions.

In the second part, we look more closely at the ``randomness'' case and further divide the proof into two subcases, depending on whether the ``typical density'' of $F_Y$ is high or low. This is a technical strengthening of a natural generalization of Balogh and Luo's idea, addressing challenges that arise only in higher dimensions. More precisely, the ``randomness'' case guarantees that for some $1\le i\le d-1$ there are $\Omega(n^{i+1})$ IGP $(i+1)$-sets $I$ such that $F_I$ is $(i,j,\varepsilon, X)$-balanced, i.e., 
\[
  |X\cap F_I| \;\in\; \bigl(\varepsilon^{\,d-i}n q^{-(j+1)/d},\ \varepsilon^{\,d-i}n q^{-j/d}\bigr]
  \quad\text{and}\quad
  |X\cap F'|\;\le\; \varepsilon\,|X\cap F_I|
\]
for every $(i-1)$-flat $F'\subseteq F_I$.

\medskip

\noindent\textbf{High density ($j\le i$).}
In this case, we generate our coplanar $(d+1)$-sets as follows. 

\medskip

\textbf{Step 1.} Pick a typical IGP $(i+1)$-set $I$ for which $F_I$ is $(i,j,\varepsilon, X)$-balanced.

\textbf{Step 2.} Pick a point $u\in F_I$ such that $I\cup\{u\}$ is critical coplanar. Since $F_I$ is pseudorandom and has high density, there are sufficiently many ways to choose $u$. 

\textbf{Step 3.} Choose arbitrarily the remaining $d-i-1$ points $W=\{w_1,\dots,w_{d-i-1}\}$. 

\textbf{Step 4.} Retain the $(d+1)$-set $I\cup\{u\}\cup W$ independently with a prescribed probability.

\medskip

\noindent The fact that $I\cup\{u\}$ is critical coplanar, together with the many choices for  $u$, guarantees that the generated collection of coplanar $(d+1)$-sets has the desired properties.

\medskip

\noindent\textbf{Low density ($j>i$).}
The real challenge lies in the ``low density'' case when $j>i$. Balogh and Luo~\cite{LB} addressed the case $d=3, j=2, i=1$ with an ingenious ad hoc method, though their approach appears difficult to extend to higher uniformities. Our new idea is to somehow reduce the ``low density'' case to the already-solved ``high density'' case. To this end, we make use of the hypergraph supersaturation lemma to ``lift'' the given $\Omega(n^{i+1})$ IGP $(i+1)$-sets into $\Omega(n^{j+1})$ IGP $(j+1)$-sets. For example, in the first non-trivial case where $(d, j, i) = (4, 3, 1)$, by definition there are $\Omega(n^2)$ pairs $\{x,y\}$ such that $|\ell_{xy}\cap X|\in (\varepsilon^3nq^{-1},\varepsilon^3nq^{-\frac{3}{4}}]$. 

\begin{figure}[H]
    \centering
    \includegraphics[width=0.6\linewidth]{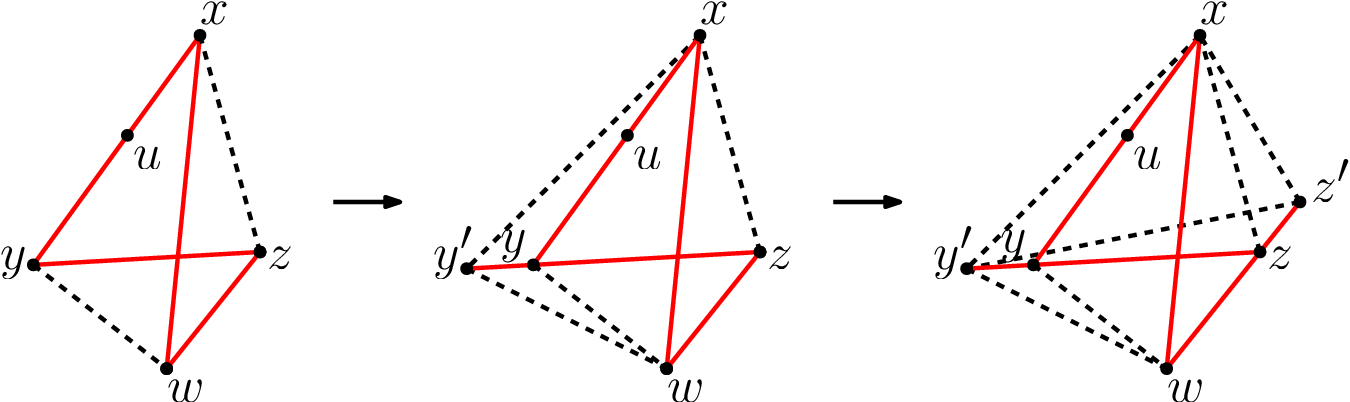}
    \caption{$\ell_{xy}$, $\ell_{yz}$, $\ell_{zw}$, $\ell_{wx}$ are the dense lines. If $u$ lies in $\ell_{xy}$, then we substitute point $y$ with some new point $y'\in \ell_{yz}$ and $z$ with some new point $z'\in \ell_{zw}$.}
    \label{fig3}
\end{figure}

We construct an auxiliary graph on $X$ whose edges are the $\Omega(n^2)$ pairs $\{x,y\}$. Since the 4-cycle $C_4$ is bipartite,  supersaturation yields $\Omega(n^4)$ 4-sets that contains a copy of $C_4$ in the auxiliary graph. We can further assume that these $\Omega(n^4)$ 4-sets $\{x,y,z,w\}$ are IGP, since by randomness the number of coplanar $4$-sets is $o(n^4)$. Moreover, the density of $F_{xyzw}$ falls into the ``high density'' realm. Next we want to pick a point $u\in F_{xyzw}$ such that $\{x,y,z,w,u\}$ is critical coplanar. However, since $F_{xyzw}$ is not ``balanced'', there might not be sufficiently many choices for such a point $u$. Suppose we randomly pick a point $u$ in $F_{xyzw}$, it is likely that $u$ lies on some dense line, say $\ell_{xy}$. In such ``bad'' events, we find substitutions for the points $x,y,z,w$ along the dense lines $\ell_{xy},\ell_{yz}, \ell_{zw},\ell_{wx}$ to obtain a critical coplanar 5-set (see Figure~\ref{fig3} for an illustration). 
Note that the support of $u$ in terms of the other four points has been enlarged during the process $\{x,y\}\rightarrow\{x,y',z\}\rightarrow\{x,y',z',w\}$, hence we call this operation the \textit{support-enlarging substitution}. Since the lines are dense (and for larger $i$, the $i$-flats are dense and balanced), such substitutions can always be found. A problem of this method is that each critical coplanar $(j+2)$-set could be generated many times, i.e., there might be overcounting. To overcome this, we construct, for each dense and balanced $i$-flat $F$, a \textit{substitution bipartite graph} $\Gamma_F$ whose two parts are both $X\cap F$ and whose maximum degree is bounded by a constant independent of $q$. During the support-enlarging substitution, if a point $x\in F$ need a substitution $x'\in F$, then we will find $x'$ in $N_{\Gamma_F}(x)$, the neighborhood of $x$ in $\Gamma_F$. Since each $x'\in F$ is adjacent to at most constantly many $x\in F$, and the number of points being substituted is at most $j$, we ensure that each critical coplanar $(j+1)$-set is overcounted only $O(1)$ times, which is acceptable. 
After constructing the critical coplanar $(j+2)$-set $I\cup\{u\}$, we repeat the same process as in the high-density case: pick the remaining $d-j-1$ points arbitrarily and retain each $(d+1)$-set independently with an appropriate probability. 
This yields the desired collection of coplanar $(d+1)$-sets, completing the proof.

\section{The structure vs. randomness lemma}\label{blanceprop}

The following lemma shows that either there are quite a lot of coplanar $(d+1)$-sets or there are many IGP $d$-sets that determine a ``dense'' flat.

\begin{lemma}\label{lemma:largeflat}
For any integer $d \ge 2$, there exists a positive constant $\varepsilon > 0$ such that the following holds for all sufficiently large prime powers $q$.  
If $X \subseteq \mathbb{F}_q^d$ is a point set with $|X|=n$, then either there are $\Omega(n^{d+1})$ coplanar $(d+1)$-sets in $X$, or there exist at least $\tfrac{1}{2}\binom{n}{d}$ IGP $d$-sets $Y$ with
\[
   |X \cap F_Y| \;\ge\; \varepsilon n q^{-1}.
\]
\end{lemma}

\begin{proof}
     Let $1/q\ll \varepsilon\ll 1/d$. Let $G_d$ be the set of IGP $(d+1)$-sets and $S_d$ be the set of coplanar $(d+1)$-sets. Then $|G_d|+|S_d|=\binom{n}{d+1}\ge (d+1)^{-d-1}n^{d+1}$. 
    
    Suppose that $|S_d|\le \varepsilon^{2d} n^{d+1}$. Then each $s$-flat with $s< d$ contains at most $\varepsilon n$ points; otherwise such a flat would contain at least $\binom{\varepsilon n}{d+1}>\varepsilon^{2d}n^{d+1}$ coplanar $(d+1)$-sets. 

    Let $S^\ast$ be the set of IGP $d$-sets $\{x_1,\dots, x_d\}$ with 
    $$|F_{\{x_1,\dots, x_d\}}\cap X|\;\ge\; \varepsilon n q^{-1}.$$
    
   By choosing \(x_1,\dots,x_{d-1}\) iteratively with \(x_i \notin F_{\{x_1,\dots, x_{i-1}\}}\) for \(i=2,\dots,d-1\), the number of IGP \((d-1)\)-sets in \(X^{d-1}\) is at least 
    \(\big((1-\varepsilon)n\big)^{d-1}/(d-1)!\). 
    Fix an IGP \(\{x_1,\dots,x_{d-1}\}\) and let $F = F_{\{x_1,\dots, x_{d-1}\}}$. 
    There are \(q+1\) distinct \((d-1)\)-flats containing \(F\), say \(F_1,\dots,F_{q+1}\).  
    Since \(|F \cap X| \le \varepsilon n\) and \(\bigcup_{i \in [q+1]} (F_i \cap X) = X\) with \(F_i \cap F_j = F\) for all \(i \ne j\), we have
    \[
        \sum_{i \in [q+1]} |(F_i \setminus F) \cap X|
            \;=\; |X| \;-\; |F| \;\ge\; (1-\varepsilon)n.
    \]
    Therefore 
 \[
        \sum_{\substack{i \in [q+1] \\ |F_i \setminus F| \,\ge\, \varepsilon n/q}}
            |(F_i \setminus F)\cap X|
            \;\ge\; (1-\varepsilon)n - (q+1)\cdot \varepsilon n/q
            \;\ge\; (1-3\varepsilon)n.
    \]
    Hence the number of $x_d$ such that $\{x_1,\dots, x_{d-1},x_d\}\in S^{\ast}$ is at least $(1-3\varepsilon)n$. Thus, \[
        |S^\ast| \;\ge\; \frac{(1-\varepsilon)^{d-1} n^{d-1}}{(d-1)!} \cdot \frac{(1-3\varepsilon)n}{d}
            \;\ge\; \frac{(1-2d\varepsilon)n^d}{d!}
            \;\ge\; \dfrac{1}{2}\binom{n}{d}.
    \]
    This completes the proof.
\end{proof}

Before proving the key lemma, we first record a structural observation that allows us to classify the distribution of points in a flat.

\begin{lemma}[Structure vs. Randomness Lemma]\label{concentrationprop}
There exist positive constants $C,\varepsilon$ such that the following holds for all sufficiently large prime powers $q$. 
If $X \subseteq \mathbb{F}_q^d$ is a point set with $|X|=n \ge Cq$ and $d \ge 2$, then either there are $\Omega(n^{d+1})$ coplanar $\l(d+1\r)$-sets in $X$, or there exist $1\le i\le d-1$ and $0\le j\le d-1$ such that there are at least $\Omega(n^{i+1})$ $(i+1)$-sets $I \subseteq X$ have the flat $F_I$ $(i,j,\varepsilon,X)$-balanced. 

Moreover, if $j> i$, then there are $\Omega(n^{j+1})$ $(i,\varepsilon,X)$-good $(j+1)$-sets in $X$. 
\end{lemma}

\begin{proof}
    Let $1/q\ll \varepsilon\ll 1/d$.  
We choose constants $\gamma$, $\varepsilon_1,\dots, \varepsilon_d$ and $\beta_0,\beta_1,\dots, \beta_{d-1}$ such that

$$1/q\ll \gamma\ll \beta_0\ll \dots \ll \beta_{d-1}\ll \varepsilon_1\ll \dots \ll \varepsilon_d\ll \varepsilon\ll 1/C.$$ 

Let $S_d$ be the set of coplanar $(d+1)$-sets and let $S^\ast$ be the set of IGP $d$-sets $\{x_1,\dots, x_{d-1},x_d\}$ satisfying
    \[
    |F_{x_1\dots x_{d-1}x_d}\cap X|\;\ge\; \varepsilon n q^{-1}.
\] By Lemma~\ref{lemma:largeflat}, we may assume that \[
    |S_d|\;\le\; \gamma^{2d} n^{d+1}
    \quad\text{and}\quad
    |S^\ast|\;\ge\; \dfrac{1}{2}\binom{n}{d}.
\] This implies that every $s$-flat with $s<d$ contains at most $\gamma n$ points; otherwise such an $s$-flat would contribute at least $\binom{\gamma n}{d+1}>\gamma^{2d} n^{d+1}$ coplanar $(d+1)$-sets, a contradiction.
    
    Now, for integers $1 \le i \le d-1$ and $0 \le j \le d-1$, let $A_{i,j}$ denote the set of IGP $(i+1)$-sets $\{x_1,\dots,x_i,x_{i+1}\}$ in $X$ such that
\[
    |F_{x_1\dots x_ix_{i+1}}\cap X|
    \;\in\; \bigl(\,\varepsilon^{\,d-i} n q^{-\frac{j+1}{d}},\;
                     \varepsilon^{\,d-i} n q^{-\frac{j}{d}}\,\bigr].
\]

Note that for each $d$-set $\{x_1,\dots, x_{d-1}, x_d\}$ in $S^\ast$, we have $\varepsilon n\ge \gamma n \ge |F_{x_1\dots x_{d-1} x_d}|> \varepsilon n/q$. This implies that \cmt{$$(x_1,\dots, x_{d-1}, x_d)\in \bigcup_{j\in [d-1]\cup \{0\}}A_{d-1,j}.$$}
$$S^*=\bigsqcup_{0\le j\le d-1}A_{d-1,j}.$$  
Therefore,
\[
   \sum_{\substack{0 \le j \le d-1\\ 1\le i\le d-1}} |A_{i,j}|\,n^{\,d-i-1}
   \;\ge\; \sum_{0\le j\le d-1}|A_{d-1,j}|
   \;=\; |S^\ast|
   \;\ge\; \dfrac{1}{2}\binom{n}{d}.
\]
By the pigeonhole principle, there exists at least one $j$ such that \[
   \sum_{1\le i\le d-1}|A_{i,j}|\,n^{\,d-i-1}
   \;\ge\; \beta_j n^d.
\] Let~$j^\ast$ be the minimum such index. Again, by the pigeonhole principle, there exists at least one $i$ such that 
$$|A_{i,j^*}|\;\ge\; \beta_{j^*}\varepsilon_i  n^{i+1}.$$ 
Let~$i^\ast$ be the minimum such index. 

\begin{claim}
    There are $\Omega(n^{i^\ast+1})$ $(i^\ast+1)$-sets~$I$ in $X$ such that the flat $F_I$ is $(i^\ast,j^\ast ,\varepsilon,X)$-balanced. 
\end{claim}
\begin{proofclaim}
    
For each $s\in [i^\ast]$, let $A_{i^\ast,j^\ast,s}$ denote the family of $(i^\ast+1)$-sets $\{x_1,\dots, x_{i^\ast},x_{i^\ast+1}\}$ in $A_{i^\ast,j^\ast}$ such that there is an $s$-flat in $F_{x_1\dots x_{i^\ast}x_{i^\ast+1}}$ with more than $\varepsilon^{i^\ast-s} |F_{x_1\dots x_{i^\ast}x_{i^\ast+1}}\cap X|$ points in $X$ and there is no $(s-1)$-flat in $F_{x_1\dots x_{i^\ast} x_{i^\ast+1}}$ with more than $\varepsilon^{i^\ast-(s-1)} |F_{x_1\dots x_{i^\ast} x_{i^\ast+1}}\cap X|$ points in $X$. Intuitively, the sets in $A_{i^\ast,j^\ast,s}$ are those $(i^\ast+1)$-sets in $A_{i^\ast,j^\ast}$ whose span has
``concentration'' at dimension $s$, but not at any lower dimension.  

We now focus on these families. By considering the number of $(s+1)$-sets lying on such
$s$-flats, we will deduce either a contradiction with the minimality of $i^\ast$ or with
that of $j^\ast$.

 Note that $A_{i^\ast,j^\ast,i^\ast}$ be the subfamily of $A_{i^\ast,j^\ast}$ consisting of $(i^\ast\!+\!1)$-sets whose span is an $(i^\ast,j^\ast,\varepsilon,X)$-balanced $i^\ast$-flat. We claim
\[
A_{i^\ast,j^\ast}=\bigcup_{s\in[i^\ast]}A_{i^\ast,j^\ast,s}.
\]
Otherwise, take $I\in A_{i^\ast,j^\ast}$ with $I\notin A_{i^\ast,j^\ast,s}$ for all $s$. In particular $I\notin A_{i^\ast,j^\ast,i^\ast}$, so there is an $(i^\ast\!-\!1)$-flat $P_{i^\ast-1}\subset F_I$ with $|P_{i^\ast-1}\cap X|>\varepsilon|F_I\cap X|$. Iterating down the dimension yields a $0$-flat $P_0$ with
\[
|P_0\cap X|>\varepsilon^{\,i^\ast}|F_I\cap X|\;\ge\;\varepsilon^{2d}\frac{n}{q}>1,
\]
a contradiction since a $0$-flat contains at most one point. Hence $A_{i^\ast,j^\ast}=\bigcup_{s\in[i^\ast]}A_{i^\ast,j^\ast,s}$. By pigeonhole,
\[
|A_{i^\ast,j^\ast,s}|\;\ge\;\frac{|A_{i^\ast,j^\ast}|}{i^\ast}\;\ge\;\frac{\beta_{j^\ast}\varepsilon_{i^\ast}}{d}\,n^{\,i^\ast+1}.
\]
If $s=i^\ast$ we are done; otherwise assume $s<i^\ast$.

    Let $\mathcal{F}$ be the set of $i^\ast$-flats determined by the $(i^\ast+1)$-sets in $A_{i^\ast,j^\ast, s}$. 
    For each $i^\ast$-flat $F$ in $\mathcal{F}$, choose an $s$-flat $P_F \subseteq F$ such that  $|X \cap P_F| \ge \varepsilon^{i^\ast-s}|F\cap X|$. 
    
    To count IGP $s$-sets lying on $P_F$ over all $F$ in $\mathcal{F}$, we construct a sequence of auxiliary graphs. 
    For $F\in \mathcal{F}$, by Lemma~\ref{auxiliary graph} let $H_F$ be a bi-$(3d^2,6d^2|F\cap X|/|P_F\cap X|)$ graph with parts $F\cap X$ and $P_F\cap X$ such that for each $v\in F\cap X$, $N(v)$ cannot be covered by $2d$ distinct $(s-1)$-flats of $P_F$.

    Now we define a map $f:A_{i^\ast,j^\ast,s}\rightarrow A_{i^\ast,j^\ast,s}$. By the function~$f$, we will map each IGP $(i^\ast+1)$-set $Y$ to a new $(i^\ast+1)$-set $Y'$ with exactly $(s+1)$ vertices on the heavy flat $P_{F_Y}$.

    Let $Y=\{x_1,\dots,x_{i^\ast+1}\}\in A_{i^\ast,j^\ast,s}$, and let $P_{F_Y}$ be an $s$-flat inside $F_Y$ with more than $\varepsilon^{\,i^\ast-s}|F_Y\cap X|$ points of $X$. Assume $P_{F_Y}$ is supported by a subset $C\subseteq Y$. Since $\dim P_{F_Y}=s$, we have $|C|=s+1$.

Define $f(Y)$ by an iterative substitution: set $Y_0=Y$ and $C_0=C$. Suppose $Y_0,\dots,Y_{t-1}$ and $C_0,\dots,C_{t-1}$ are defined. If $C_{t-1}\subseteq P_{F_Y}$, set $f(Y)=Y_{t-1}$. Otherwise pick $x\in C_{t-1}\setminus P_{F_Y}$ and choose $x'\in P_{F_Y}\cap X$ so that $x'\in N_{H_{F_Y}}(x)$ and $x'\notin F_{Y_{t-1}\setminus\{x\}}$. Let
\[
Y_t=(Y_{t-1}\cup\{x'\})\setminus\{x\}.
\]
By Lemma~\ref{propflat}\,\ref{propflat2}, each $Y_t$ is an IGP $(i^\ast\!+\!1)$-set with $F_{Y_t}=F_{Y_{t-1}}=F_Y$, and thus $P_{F_Y}\subseteq F_{Y_t}$. Let $C_t$ be the support of $P_{F_Y}$ inside $Y_t$.

\begin{claim}
    The iteration satisfies:
\begin{enumerate}[label=\rm{(\roman*)}]
    \item\label{claimpro1} A valid $x'$ always exists.
    \item\label{claimpro2} The process stops in at most $s+1$ steps.
    \item\label{claimpro3} For each $Z\in f(A_{i^\ast,j^\ast,s})$, $|f^{-1}(Z)|\le (4d)^{3s}\varepsilon^{(s-i^\ast)(s+1)}$.
\end{enumerate}
\end{claim}
\begin{proofclaim}
To see~\ref{claimpro1}, the existence of $x'$, by the construction of $H_{F_Y}$, it suffices to prove that
 $\dim(F_{Y_{t-1}\setminus \{x\}}\cap P_{F_Y})\le s-1$. By Lemma~\ref{propflat}~\ref{propflat3}, it suffices to show that $P_{F_Y}\nsubseteq F_{Y_{t-1}\setminus \{x\}}$. Suppose not, then $P_{F_Y}$ is supported by a subset $C'$ of $Y_{t-1}\setminus\{x\}$. This implies that $P_{F_Y}$ is supported by distinct sets $C'$ and $C_t$. This contradicts Lemma~\ref{exactprop}~\ref{exactprop2}, which asserts that the support of $P_{F_Y}$ in $Y$ is unique. Thus, $x'$ always exists.  

To see~\ref{claimpro2}, note that $P_{F_{Y}}\cap Y_t\subseteq C_t$. Otherwise there exists some point $y\in P_{F_Y}\cap Y_t$ not in $C_t$. Then by the definition of support, $y\in \text{span}(C_t)\subseteq\sp{Y_t\setminus \{y\}}$. This contradicts the fact that $Y_t$ is IGP. In the process, Each substitution increases $|P_{F_{Y}}\cap Y_t|$ by $1$ and $Y_t$ remains IGP, so $|P_{F_Y}\cap Y_t|\le \dim(P_{F_Y})+1=s+1$. The process ends once $|P_{F_Y}\cap Y_t|=s+1$, hence in at most $s+1$ steps. 

To see~\ref{claimpro3}, fix $Z$ and reverse the process. Each newly added $v\in Z\cap P_{F_Y}$ could have come from at most $6d^2\,{|F_Y\cap X|}/{|P_{F_Y}\cap X|}$ preimages. Since $|P_{F_Y}\cap X|/|F_Y\cap X|\ge \varepsilon^{\,i^\ast-s}$ by the definition of $A_{i^*, j^*, s}$, after at most $s+1$ steps we get 
\[f^{-1}(Z) \;\le\; (6d^2)^{s+1}\varepsilon^{-(i^\ast-s)(s+1)} \;\le\; (4d)^{3s}\varepsilon^{(s-i^\ast)(s+1)}.\]
\end{proofclaim}
Let \[
B_{s}:=\{f(Y)\cap P_{F_Y}:\ Y\in A_{i^*,j^*,s}\}.\]

Note that by definition every element of $B_s$ is an $(s+1)$-set and any such set is contained in at most $n^{i^*-s}$ $(i^*+1)$-sets. Then 
$$|B_s|\;\ge\; \frac{|f(A_{i^\ast,j^\ast, s})|}{n^{i^\ast-s}}\;\ge\; \frac{|A_{i^\ast, j^\ast, s}|}{(4d)^{3s}\varepsilon^{(s-i^\ast)(s+1)}n^{i^*-s}}\;\ge\; 2\varepsilon_s\,\beta_{j^\ast}\, n^{s+1}.$$
For each $S \in B_s$, $S$ determines an $s$-flat $F_S$ and 
$$
\varepsilon^{i^\ast-s}\varepsilon^{d-i^\ast} n q^{-\frac{j^\ast+1}{d}}\;\le \;|F_S|\;\le\; \varepsilon^{d-i^\ast} nq^{-\frac{j^\ast}{d}}.
$$
Define
\[
\begin{aligned}
B_s^\ast&:=\bigl\{S\in B_s:\ |F_S\cap X|\in\bigl(\varepsilon^{\,d-s}nq^{-\frac{j^\ast+1}{d}},\,\varepsilon^{\,d-s}nq^{-\frac{j^\ast}{d}}\bigr]\bigr\},\\
B_s^{\ast\ast}&:=\bigl\{S\in B_s:\ |F_S\cap X|\in\bigl(\varepsilon^{\,d-s}nq^{-\frac{j^\ast}{d}},\,\varepsilon^{\,d-s}nq^{-\frac{j^\ast-1}{d}}\bigr]\bigr\}.
\end{aligned}
\]
Since $1/q\ll\varepsilon$, we have $B_s=B_s^\ast\cup B_s^{\ast\ast}$. Thus either
\[
|B_s^\ast|\ge \varepsilon_s\beta_{j^\ast}n^{\,s+1}\quad\text{or}\quad
|B_s^{\ast\ast}|\ge \varepsilon_s\beta_{j^\ast}n^{\,s+1}.
\]
In the first case $|A_{s,j^\ast}|\ge |B_s^\ast|\ge \varepsilon_s\beta_{j^\ast}n^{\,s+1}$, contradicting the minimality of $i^\ast$.
In the second case,
\[
\sum_{i\in[d-1]}|A_{i,j^\ast-1}|\,n^{\,d-i-1}\ \ge\ |A_{s,j^\ast-1}|\,n^{\,d-s-1}
\ \ge\ |B_s^{\ast\ast}|\,n^{\,d-s-1}
\ \ge\ \varepsilon_s\beta_{j^\ast}\,n^{\,d}\gg \beta_{j^\ast-1}\,n^{\,d},
\]
contradicting the minimality of $j^\ast$. Hence $s=i^\ast$, and there are $\Omega(n^{\,i^\ast+1})$ $(i^\ast\!+\!1)$-sets $I$ with $F_I$ $(i^\ast,j^\ast,\varepsilon,X)$-balanced.
\end{proofclaim}

To show the moreover statement, for $d-1\ge j>i\ge 1$, let $J_{i,j}$ denote the complete $(i+1)$-partite $(i+1)$-graph whose part sizes are $\lceil (j-i)/2\rceil+1,\ \lfloor (j-i)/2\rfloor+1,\,1,\dots,1$.
Let $H_{i^\ast+1}(X)$ be the $(i^\ast+1)$-graph on $X$ whose edge set consists of all $(i^\ast,j^\ast,\varepsilon,X)$-sets in $X$.
By definition,
\[
|E(H_{i^\ast+1})|\;\ge\;\beta_{j^\ast}\varepsilon_{i^\ast}\,n^{\,i^\ast+1}.
\]
Applying Lemma~\ref{supersaturation} with $F=J_{i^\ast,j^\ast}$ and $\varepsilon=\beta_{j^\ast}\varepsilon_{i^\ast}$, there exist at least $\delta n^{j^\ast+1}$ copies of $J_{i^\ast,j^\ast}$, where $\delta=\delta(\varepsilon_{i^\ast},\beta_{j^\ast}) > 0$. At least $\delta n^{\,j^\ast+1}/2$ of these copies have IGP vertex sets; otherwise, at least $\delta n^{\,j^\ast+1}/2$ of them are coplanar $(j^\ast+1)$-sets, which then yield
\[
\Omega\bigl((\delta/2)\,n^{\,j^\ast+1}\cdot n^{\,d-j^\ast}\bigr)
=\Omega(\delta\,n^{\,d+1})
\]
coplanar $(d+1)$-sets, contradicting $ |S_d|\le \gamma n^{\,d+1}$ since $\gamma\ll \varepsilon_{i^\ast},\beta_{j^\ast},1/d$.

Note that at least $\delta n^{\,j^\ast+1}/4$ copies of $J_{i^\ast,j^\ast}$ correspond to $(j^\ast+1)$-sets that are $(i^\ast,\varepsilon,X)$-good. Indeed, if fewer than that many were $(i^\ast,\varepsilon,X)$-good, then at least $\delta n^{\,j^\ast+1}/4$ of the IGP $(j^\ast+1)$-sets $J$ would satisfy
\[
n\,q^{-j^\ast/d}\ \le\ |F_J\cap X|\ \le\ \gamma n\ \ll\ \varepsilon^{\,d-j^\ast}n,
\]
and hence $J\in \bigcup_{t<j^\ast}A_{j^\ast,t}$.
By the pigeonhole principle, there exists $t<j^\ast$ with
\[
|A_{j^\ast,t}|\ \ge\ \frac{\delta}{4d}\,n^{\,j^\ast+1}.
\]
Consequently,
\[
\sum_{i=1}^{d-1}|A_{i,t}|\,n^{\,d-i-1}
\ \ge\ |A_{j^\ast,t}|\,n^{\,d-j^\ast-1}
\ \ge\ \frac{\delta}{4d}\,n^{\,j^\ast+1}\,n^{\,d-j^\ast-1}
\ =\ \frac{\delta}{4d}\,n^{\,d}
\ \gg\ \beta_t\,n^{\,d},
\]
contradicting the minimality of $j^\ast$.
Therefore there are $\Omega(n^{\,j^\ast+1})$ $(i^\ast,\varepsilon,X)$-good $(j^\ast+1)$-sets, as claimed.
\end{proof}

\section{The ``structure'' case}\label{sec:struc}
In this section, we prove the structure case, namely, the situation in which $X$ contains $\Omega(n^{d+1})$ coplanar $(d+1)$-sets. 
\begin{lemma}\label{struc}
    There exists a constant $C$ such that the following statement holds. Let $X$ be a point set of $\mathbb{F}_q^d$ with $|X|=n\ge Cq$. Suppose that there are $\Omega(n^{d+1})$ coplanar $(d+1)$-sets in $X$. Then there exists a collection $\mathcal{S}$ of coplanar $(d+1)$-sets in $X$ such that 
     \begin{enumerate}[label=\rm{(\roman*)}]
        \item\label{equation11} $|\mathcal{S}|=\Omega\l(n^{d+1}/q\r)$ and
        \item\label{equation21} $\Delta_s(\mathcal{S})\le O\l(n^{d+1-s}\cdot{q^{-\frac{d+1-s}{d}}}\r)$ for $s\in [d+1]$. 
    \end{enumerate}
\end{lemma}
\begin{proof}
    Let $\mathcal{C}$ be the collection of all coplanar $(d+1)$-sets in $X$. Select each $C\in \mathcal{C}$ with probability $p=q^{-1}$ and let $\mathcal{S}$ be the resulting random collection.

 Since $|\mathcal{C}| = \Omega(n^{d+1})$, we have
    \[
\mathbb{E}[|\mathcal{S}|] 
   \;\ge\; \frac{|\mathcal{C}|}{q} 
   \;=\; \Omega\!\left( \frac{n^{\,d+1}}{q} \right).
\]
    By the Chernoff bound (Lemma~\ref{lem: Chernoff}), $$\Pr\left[\,|\mathcal S|\le \mathbb{E}[|\mathcal S|]/2\,\right] \;\le\; \exp(-\Omega(\mathbb{E}[\mathcal{S}]))\;\le\; o(1).$$
    Hence, with probability $1-o(1)$, we have
    \[|\S|\;\ge\;\Omega\left(\frac{n^{d+1}}{q}\right),\]
    which establishes property \ref{equation11}.

    Now, we check~\ref{equation21}. We first bound expectations. We will show that for any $d$-set $D$ in $X$,
    \begin{align}\label{Delta41}
         \E[\deg(D)] \;\le\; O\l(n \cdot q^{-\frac{d-1}{d}}\r),
    \end{align}
    and for any $1$-set $D$ in $X$
    \begin{align}\label{Delta42}
        \E[\deg(D)]\;\le\; O\l(n^{d}\cdot q^{-1}\r).
    \end{align}
     Note that, for any $1\le s\le s'\le d$, any $s$-set $D$ satisfies 
     \[\deg(D)\;\le\; n^{s'-s}\max\l\{\deg(D'):\ D\subset D',\,|D'|=s\,\r\}.\] 

     Thus~(\ref{Delta41}) and~(\ref{Delta42}) will imply, for any $1\le s\le d$ and any $s$-set $D$,
     \begin{align}\label{DeltaExpectation0}
     \E[\deg(D)]\;\le\; O\l(n^{d+1-s}\,q^{-\frac{d+1-s}{d}}\r).
     \end{align}

    To see~(\ref{Delta41}), we fixed a $d$-set $D$. Then there are at most $n$ potential choices $v$ such that $D\cup \{v\}\in \mathcal{C}$. Since we pick each set in $\mathcal{C}$ with probability $1/q$. This implies that $\mathbb{E}[\deg(D)]\le \frac{n}{q}$ as desired. 
    
     {To see~(\ref{Delta42}), we fixed an $1$-set $D$. Then there are at most $n^{d}$ potential choices $C$ such that $D\cup C\in \mathcal{C}$. Since we pick each set in $\mathcal{C}$ with probability $1/q$. This implies that $\mathbb{E}[\deg(D)]\le n^d/ q$ as desired.}
     
     Now we have inequalities~(\ref{DeltaExpectation0}). Thus for any $1\le s\le d$, there exists a sufficiently large constant $C_s$ such that for any $s$-set $D$ in $X$, we have
    $\E[\deg(D)]\le C_sn^{d+1-s}q^{-\frac{d+1-s}{d}}.$ For any fixed $s$-set $D$, note that $\deg(D)$ is a sum of independent Bernoulli random variables. Thus by the Chernoff bound (Lemma~\ref{lem: Chernoff}), together with the fact that $q\le O(n)$, $d\ge 2$ and $s\le d$, there exists a constant $C_s$ such that 
    $$
    \Pr\l[\deg(D)\;\ge\; 20C_s\, n^{d+1-s}\,q^{-\frac{d+1-s}{d}}\r]\;\le\; 2\exp\l(-10C_s\, n^{d+1-s}\,q^{-\frac{d+1-s}{d}}\r)\;\le\; \exp\l(-\Omega\l(n^{\frac{d-1}{d}}\r)\r).
    $$

    By the union bound, for any $1\le s\le d$,
    $$
    \Pr\l[\Delta_{s}(\mathcal{S})\;\ge\; 20C_s\, n^{d+1-s}\,q^{-\frac{d+1-s}{d}}\r]\;\le\; \binom{n}{s}\exp\l(-\Omega\l(n^{\frac{d-1}{d}}\r)\r)\;\le\; o(1).
    $$
     
     Hence, with probability $1-o(1)$, we have $\Delta_s(\mathcal{S})\le O\l(n^{d+1-s}q^{-\frac{d+1-s}{d}}\r)$.

    Therefore, we can find the collection $\mathcal{S}$ as desired.
\end{proof}

\section{The ``randomness'' case with high density}\label{BH}
In this section, we prove the ``high density'' subcase in the ``randomness'' case, where there are $\Omega(n^{i+1})$ IGP $(i+1)$-sets determining a balanced flat with high density.
\begin{lemma}\label{high}
    There exists a constant $C$ such that the following statement holds. Let $X$ be a point set of $\mathbb{F}_q^d$ with $|X|=n\ge Cq$. Suppose that there are $o(n^{d+1})$ coplanar $(d+1)$-sets in $X$ and for some $i\in [d-1], 0\le j\le i$ and $1/q\ll \varepsilon\ll 1/d$, there are at least $\Omega(n^{i+1})$ $(i+1)$-sets $I$ in $X$  such that the flat $F_I$ is $(i,j,\varepsilon,X)$-balanced. Then there exists a collection $\mathcal{S}$ of coplanar $(d+1)$-sets in $X$ such that 
     \begin{enumerate}[label=\rm{(\roman*)}]
        \item\label{equation1} $|\mathcal{S}|=\Omega\l(n^{d+1}/q\r)$ and
        \item\label{equation2} $\Delta_s(\mathcal{S})\le O\l(n^{d+1-s}\cdot{q^{-\frac{d+1-s}{d}}}\r)$ for $s\in [d+1]$. 
    \end{enumerate}
\end{lemma}
\begin{proof}
Let $\mathcal{I}$ be the family of  $(i+1)$-sets $I$ in $X$ with $F_I$ $(i,j,\varepsilon,X)$-balanced. Then 
\[|F_I\cap X|\;\ge\; \varepsilon^{d-i}\,n\,q^{-\frac{j+1}{d}}\;\ge\; \varepsilon^{d-i}\,n\,q^{-\frac{i+1}{d}}.\]

We generate $\mathcal{S}$ as follows.

\medskip

\textbf{Step 1.} Arbitrarily choose an IGP $(i+1)$-set $I=\{x_1,\dots,x_{i+1}\}$ from $\mathcal{I}$. Let the distinct $(i-1)$-flats determined by $i$-subsets of $I$ be $F_1,\dots, F_{i+1}$, and write $F := F_I$. Then $|F\cap X|\ge \varepsilon^{d-i}nq^{-\frac{i+1}{d}}$. 

\textbf{Step 2.} Arbitrarily choose a point $u$ from $\l(F\setminus\bigcup_{s\in [i+1]} F_s\r)\cap X$. Let $U=I\cup \{u\}$.

\textbf{Step 3.} Arbitrarily choose the remaining $d-i-1$ points $W=\{w_1,\dots, w_{d-i-1}\}\subset X\setminus U$.

\textbf{Step 4.} Keep the coplanar $(d+1)$-set $U\cup W$ with probability $\dfrac{\varepsilon^{d-i} n}{q|F\cap X|}$.

\medskip

Note that since $F$ is $(i,\varepsilon,X)$-balanced, the number of choices for $u\in (F\setminus \bigcup_{s\in [i+1]}F_s)\cap X$ is at least
$$
(1-(i+1)\varepsilon)\ |F\cap X|\;\ge\; \frac{1}{2}|F\cap X|.
$$ 

Hence
\[
\begin{aligned}
   \mathbb{E}[|\mathcal{S}|]
   &\;\ge\;
   |\mathcal{I}| \cdot 
   \bigl|\,(F \setminus \!\bigcup_{s \in [i+1]} F_s)\cap X\,\bigr|\,
   (n-d)^{\,d-i-1}\,
   \frac{\varepsilon^{\,d-i} n}{q\,|F \cap X|}\,
   \frac{1}{d!} \\[6pt]
   &\;\ge\;
   |\mathcal{I}| \cdot \dfrac{1}{2}\,|F \cap X|\,
   \left(\dfrac{n}{2}\right)^{d-i-1}\,
   \frac{\varepsilon^{\,d-i} n}{q\,|F \cap X|}\,
   \frac{1}{d!} \\[6pt]
   &\;\ge\;
   \Omega\!\left(\dfrac{n^{\,d+1}}{q}\right).
\end{aligned}
\]

Note that $|\mathcal{S}|$ is a sum of independent Bernoulli random variables. Thus by the Chernoff bound (\Cref{lem: Chernoff}),
$$
\Pr\l[\,|\mathcal{S}|\le \E[|\mathcal{S}|]/2\,\r]\;\le\; \exp\l(-\Omega(\E[|\mathcal{S}|])\r)\;\le\; o(1).
$$
Hence, \ref{equation1} holds with probability $1-o(1)$.

 Now we check~\ref{equation2}. We first check the expectation. To this end, we will only show that, for any $d$-set $D$ in $X$,
    \begin{align}\label{Delta4}
         \E[\deg(D)]\;\le\; O\l(nq^{-\frac{d-i-1}{d}}\r),
    \end{align}
   and, for any $(i+1)$-set $D$ in $X$,
    \begin{align}\label{Delta5}
        \E[\deg(D)]\;\le\; O\l(n^{d-i-1}q^{-1}\r).
    \end{align}
     Note that, for any $1\le s\le s'\le d$, any $s$-set $D$ satisfies    \[\deg(D)\;\le\; n^{s'-s}\max\l\{\deg(D'):\ D\subset D',\,|D'|=s\,\r\}.\]  Thus~(\ref{Delta4}) and (\ref{Delta5}) would imply, for any $1\le s\le d$ and any $s$-set $D$,
     \begin{align}\label{DeltaExpectation}
     \E[\deg(D)]\;\le\; O\l(n^{d+1-s}\,q^{-\frac{d+1-s}{d}}\r).
     \end{align}
     
     First, to see~(\ref{Delta4}), note that for a fixed $d$-set $D$, the remaining point have at most $n$ choices, and in step $4$, we keep the $(d+1)$-set with probability $\frac{\varepsilon^{d-i} n}{q|F\cap X|}$. Recall that $|F\cap X|\ge \varepsilon^{d-i}nq^{-\frac{i+1}{d}}$. Thus the expectation of $\deg(D)$ can be estimated as follows.
    \[
   \mathbb{E}[\deg(D)]
   \;\le\;
   n \cdot \frac{\varepsilon^{\,d-i}\,n}{\,q\,|F \cap X|}
   \;\le\;
   \frac{n}{q^{\frac{d-i-1}{d}}}.
\]

    To see~(\ref{Delta5}), we fixed a $(i+1)$-set $D$. We need to count the number of sets $U\cup W=\{x_1,\dots,x_{i+1},u,w_1,\dots,w_{d-i-1}\}$ generated by the steps above that contains $D$. First we decide which $i+1$ points in $U\cup W$ is specified by $D$ --- there are $\binom{d+1}{i+1}$ ways to do this. Note that the critical coplanar $(i+2)$ set $U=\{x_1,\dots,x_{i+1},u\}$ cannot be contained in $D$ (since $i+2>i+1$). Thus there exists a point $v\in U$ not specified. Next we specify all other $d-i-1$ unspecified points except $v$ --- there are at most $n^{d-i-1}$ choices. Since $U$ is critical coplanar, after all $i+1$ points except $v$ in $U$ have been specified, we have specified the $i$-flat $F$ which must contain $v$. So there are at most $|F\cap X|$ ways to specify $v$. Hence, we can estimate the expectation of $\deg(D)$ as follows.
    \[
   \mathbb{E}[\deg(D)]
   \;\le\;
   \binom{d+1}{\,i+1}\,
   n^{\,d-i-1}\,|F \cap X|\,
   \frac{\varepsilon^{\,d-i}\,n}{q\,|F \cap X|}
   \;=\;
   O\!\left(\dfrac{n^{\,d-i}}{q}\right).
\]
    Now we have inequalities~(\ref{DeltaExpectation}). Thus for any $1\le s\le d$, there exists a sufficiently large constant $C_s$ such that for any $s$-set $D$ in $X$, we have
    $\E[\deg(D)]\le C_sn^{d+1-s}q^{-\frac{d+1-s}{d}}.$ For any fixed $s$-set $D$, note that $\deg(D)$ is a sum of independent Bernoulli random variables. Thus by the Chernoff bound (Lemma~\ref{lem: Chernoff}), together with the fact that $q\le O(n)$, $d\ge 2$ and $s\le d$, 
    \[
    \Pr\l[\deg(D)\;\ge\; 20C_s\, n^{d+1-s}\,q^{-\frac{d+1-s}{d}}\r]\;\le\; 2\exp\l(-10C_s\, n^{d+1-s}\,q^{-\frac{d+1-s}{d}}\r)\;\le\; \exp\l(-\Omega\l(n^{\frac{d-1}{d}}\r)\r).
    \]

    By the union bound, for any $1\le s\le d$,
    $$
    \Pr\l[\Delta_{s}(\mathcal{S})\;\ge\; 20C_s\, n^{d+1-s}\,q^{-\frac{d+1-s}{d}}\r]\;\le\; \binom{n}{s}\exp\l(-\Omega\l(n^{\frac{d-1}{d}}\r)\r)\;\le\; o(1).
    $$
     
     Hence, with probability $1-o(1)$, we have $\Delta_s(\mathcal{S})\le O(n^{d+1-s}q^{-\frac{d+1-s}{d}})$.

    Therefore, we can find the collection $\mathcal{S}$ as desired.
\end{proof}

\section{The ``randomness'' case with low density}\label{BL}
In this section, we handle the ``low density'' subcase of the ``randomness'' case, where there are $\Omega(n^{i+1})$ IGP $(i+1)$-sets determining balanced but relatively ``sparse'' flats.

\subsection{Support-enlarging substitution}

We first describe a one-step substitution that enlarges support.

\begin{lemma}\label{map sublemma}
    There exists a constant $C$ such that the following statement holds. Let $X$ be a point set of $\mathbb{F}_q^d$ with  $|X|=n\ge Cq$. 
    For $1\le i<j\le d$, let $J$ be an IGP $(j+1)$-set in $X$ and let $I\subseteq J$ be an $(i+1)$-set. 
Let $y\in J$. 
Suppose $u\in F_J$ is supported by a set $U$ with $I\nsubseteq U\subseteq J$, and let $x\in U\cap I$.

    Then there exists a collection $\mathcal{I}$ of flats in $F_I$, each of dimension at most $i-1$, with $|\mathcal{I}|\le d+2$, such that for every\[
  x'\in F_I\setminus \bigcup_{F\in\mathcal{I}} F,
  \qquad
  J' := (J\setminus\{x\})\cup\{x'\},\quad
  I' := (I\setminus\{x\})\cup\{x'\},
\] we have 
    \begin{enumerate}[label=\rm{(M\arabic*)}]
        \item\label{map sublemma1} $u$ is supported by $U\cup I'$ in $J'$;
        \item\label{map sublemma2} $F_{J'}=F_{J}$;
        \item\label{map sublemma3} $F_{I'}=F_I$;
        \item\label{map sublemma4} $F_{J'\setminus \{y'\}}=F_{J\setminus \{y\}}$ for some $y'\in J'$.
    \end{enumerate}
\end{lemma}

\begin{proof}
Let $\mathcal{I}$ be the collection consisting of the flats \[
   F_I \cap F_{\,U \cup I \cup \{u\} \setminus \{z,x\}} \quad (z \in I \cup U \setminus \{x\}), 
   \qquad 
   F_{J \setminus \{x\}} \cap F_I,
   \qquad 
   F_{I \setminus \{x\}} .
\]  By~Lemma~\ref{enlargecoverset}~\ref{enlargecoverset1} and \ref{enlargecoverset2} with $(U,I,u,x,z)=(U,Y,u,w,y)$, we have
\[
\dim(F_{U\,\cup\, I\,\cup\, \{u\}\setminus\{x,z\}}\cap F_{I})\;<\;\dim(F_{I})\;=\;i,
\]
and also
$\dim(F_{J\setminus\{x\}}\cap F_I)=i-1$ and $\dim(F_{I\setminus \{x\}})=i-1$. Hence $|\mathcal{I}| \le d+2$, and every member has  dimension at most $i-1$.

Next, for each point $x'\in F_I\setminus \bigcup_{F\in \mathcal{I}}F$, by Lemma~\ref{enlargecoverset}~\ref{enlargecoverset3} with $(U,I,u,x,z)=(U,Y,u,w,y)$, $u$ is supported by $U\cup I\cup \{x'\}\setminus \{x\}=U\cup I'$. Thus~\ref{map sublemma1} holds.

For~\ref{map sublemma2}, since $x'\notin F_{J\setminus \{x\}}$ and $J$ is IGP, Lemma~\ref{propflat}~\ref{propflat2} implies that $J'$ is IGP and $F_{J'}=F_J$.

For~\ref{map sublemma3}, note that $I\subseteq J$ implies $I$ is IGP. Since $x'\notin F_{I\setminus \{x\}}$ and $I$ is IGP, again by Lemma~\ref{propflat}~\ref{propflat2}, we have $I'$ is IGP and $F_{I'}=F_I$.

To see~\ref{map sublemma4}, if $x=y$, then $J'\setminus\{x'\}=J\setminus \{x\}$. Taking $y'=x'$ proves the claim. If $x\neq y$, since $F_{J}=F_{J'}$, we have $F_{J\setminus \{y\}}=F_{J'\setminus \{y\}}$, and the claim follows. 
\end{proof}

 {Recall that given $\varepsilon$ and the point set $X$,  $\mathcal{J}_{i,j}$ is the set of $(i,\varepsilon,X)$-good $(j+1)$-sets, and $f:\mathcal{J}_{i,j}\rightarrow X$ is the root function of $\mathcal{J}_{i,j}$ which maps each $J\in \mathcal{J}_{i,j}$ to its root.  For each $u\in X$, let $\mathcal{J}_{u,i,j}$ be the collection of $J$ in $\mathcal{J}_{i,j}$ such that $u\in (F_J\setminus F_{J\setminus \{f(J)\}})\cap X$}. Note that the support of $u$ in $J$ contains~$f(J)$.

We construct $g_u: \mathcal{J}_{u,i,j} \to \binom{X}{j+1}$ such that for any $J\in \mathcal{J}_{u,i,j}$, $g_u(J)\cup \{u\}$ is a critical $(j+2)$-coplanar set. Note that if $J$ is a support of $u$, then $J\cup\{u\}$ is critical coplanar. The idea is to enlarge the support of $u$ in $J$ by iterative substitution.
\begin{lemma}\label{map lemma}
    There exists a constant $C$ such that the following statement holds. Let $X \subseteq \mathbb{F}_q^d$ with $|X|=n\ge Cq$ and $1/q\ll \varepsilon\ll 1/d$. For each $u\in X$,  there exists  $g_u:\mathcal{J}_{u,i,j}\rightarrow \binom{X}{j+1}$ such that for any $J\in \mathcal{J}_{u,i,j}$:
    \begin{enumerate}[label=\rm{(A\arabic*)}]
        \item\label{A1} $g_u(J)\cup \{u\}$ is a critical coplanar $(j+2)$-set.
        \item\label{A2} $F_{g_u(J)}=F_J$ and $F_{J\setminus \{f(J)\}}=F_{g_u(J)\setminus \{y'\}}$ for some $y'\in g_u(J)$.
        \item\label{A3} $|g_u^{-1}(\tilde{J})|\le (12d)^{5d}$ for any $\tilde{J}\in g_u(\mathcal{J}_{u,i,j})$.
    \end{enumerate}
\end{lemma}

\begin{proof}
    We construct the function by iteratively substituting elements of $J$, thereby enlarging the support of $u$ in $J$. Let $\mathcal{I}$ be the set of $(i+1)$-sets $I$ with $F_{I}$ being $(i,j,\varepsilon, X)$-balanced. Let $\mathcal{F}$ be the set of $i$-flats determined by the $(i+1)$-sets in $I$.  For each $F\in \mathcal{F}$, by~Lemma~\ref{auxiliary graph} with $s=i$ and $F'=F$, we can define an auxiliary bi-$(3d^2, 6d^2)$ graph $G_F$ with parts $A_F=F\cap X$ and $B_F=F\cap X$ such that, for each $v\in A_F$, the neighborhood $N_{B_F}(v)$ cannot be covered by $d+2$ distinct $(i-1)$-flats. 
    
    Suppose $J\in \mathcal{J}_{u,i,j}$ and $J=\left\{y,x_2,\dots, x_{j+1}\right\}$ with $y=f(J)$ and $u\in F_J\setminus F_{J\setminus \{y\}}$. Assume that $u$ is supported by $U\subseteq J$.
    Let the $(i+1)$ parts of $J$ be $V_1,\dots,V_{i+1}$ with
\[
   |V_1| \;=\; \left\lceil \tfrac{j-i}{2} \right\rceil + 1,\qquad
   |V_2| \;=\; \left\lfloor \tfrac{j-i}{2} \right\rfloor + 1,\qquad
   |V_s| \;=\; 1 \ \text{ for } 3 \le s \le i+1 .
\]
    By the definitions of $f$ and $U$, we see that $y\in V_1$ and $y\in U$. 

    First, if $U=J$, then by~Lemma~\ref{exactprop}~\ref{exactprop3}, $J\cup \{u\}$ is already a critical $(j+1)$-coplanar set. Set $g_u(J)=J$. If $U\subsetneq J$, we distinguish two cases:

    \textbf{Case 1:} If $V_1\cup V_2\subseteq U$, choose $w\in V_1\setminus \{y\}$ (well-defined since $|V_1|\ge 2$). By the definition of $J$, there exists an $(i+1)$-set $I$ with $I\nsubseteq U$ and $w \in I$, such that $F_I$ is $(i,j,\varepsilon,X)$-balanced. Since $w\in U$, by Lemma~\ref{map sublemma} with $(J,I,x,y)=(J,I,w,y)$, there exists a collection $\mathcal{I}$ of flats in $F_I$ with dimension at most $i-1$.

    By the definition of the auxiliary graph $G_{F_I}$, there exists a neighbor $w'\in (X\cap F_I)\setminus \bigcup_{F\in \mathcal{I}}F$ of $w$. Consider the new $(j+1)$-set $J':=\l(J\setminus \{w\}\r)\cup \{w'\}$. By Lemma~\ref{map sublemma}, in $J'$ the point $u$ is supported by \[\l(U\cup I\cup \{w'\}\r)\setminus \{w\}\;=\;\l(V_1\cup V_2\cup I\cup \{w'\}\r)\setminus \{w\} \;=\;\l(J\setminus \{w\}\r)\cup \{w'\}=J',\] and moreover $F_{J'}=F_{J}$ and $F_{J'\setminus \{y'\}}=F_{J\setminus \{y\}}$ for some $y'\in J'$. Set $g_u(J)=J'$, which satisfies~\ref{A1} and~\ref{A2}.

    \textbf{Case 2:} If $V_1\cup V_2\nsubseteq U$, we need to the substitution in Case~1 to enlarge the support. By~Lemma~\ref{prop:deirectedhamiltonianpath}, since $|V_2|\le |V_1|\le |V_2|+1$, there exists a directed Hamiltonian path $P$ on $V_1\cup V_2$ starting from $y$. For $s\geq1$, let $J_s$ and $U_s$ be the $(j+1)$-set and the support set after $s$ iterations respectively, with $J_0=J$ and $U_0=U$. Let $W_s$ record points removed from the support and $W'_s$ record new points added; set $W_0=W'_0=\varnothing$ and $y_0=y$.
    
        Assume $J_{s-1},U_{s-1},W_{s-1},W_{s-1}',y_{s-1}$ has been defined for some $s\ge 1$ with $U\subseteq U_{s-1}\cup W_{s-1}$. If $U_{s-1}=J_{s-1}$, then $u$ is supported by $J_{s-1}$ and we set $g_u(J)=J_{s-1}$. Otherwise, let $\overrightarrow{w_sv_s}$ be the first edge of path $P$ with $w_s\in U_{s-1}$ and $v_s\notin U_{s-1}$ (existence proved below). There exists an $(i+1)$-set $I_{s-1} \subseteq J_{s-1}\cap J$ containing $w_s,v_s$ with $F_{I_{s-1}}$ $(i,j,\varepsilon,X)$-balanced.  By Lemma~\ref{map sublemma} with $(J,I,x,y)=(J_{s-1},I_{s-1},w,y_{s-1})$, there exists a collection $\mathcal{I}$ of flats in $F_{I_{s-1}}$ of dimension at most $i-1$, with $|\mathcal{I}|\le d+2$.

    By the definition of $G_{F_{I_{s-1}}}$, there exists a neighbor $w_s'\in (X\cap F_{I_{s-1}})\setminus \bigcup_{F\in \mathcal{I}}F$ of $w$. Define $J_s:=\l(J_{s-1}\setminus \{w_s\}\r)\cup \{w_s'\}$. By Lemma~\ref{map sublemma}, in $J_s$ the point $u$ is supported by \[U_s\;:=\;\l(U_{s-1}\cup I_{s-1}\cup \{w_s'\}\r)\setminus \{w_s\}\;\subseteq\; J_s,\] and moreover $F_{J_s}=F_{J_{s-1}}$ and $F_{J_s\setminus \{y_s\}}=F_{J_{s-1}\setminus \{y_{s-1}\}}$ for some $y_s\in J_s$.

Note that $v_s \notin U_{s-1}\cup W_{s-1}$ but $v_s \in I_{s-1} \subseteq U_s$, hence
\[
   |U_s| \;=\; \bigl|(U_{s-1} \cup I_{s-1} \cup \{w'_s\}) \setminus \{w_s\}\bigr|
   \;\ge\; |U_{s-1}| + 1 \;>\; |U_{s-1}|.
\]
Since $|J_s| = j+1$ for all $s$, the process terminates in at most $j+1$ rounds when $|U_s|=j+1$, yielding $U_s=J_s$.



    We verify well-definedness. 
For $s=1$, we have $W_{0}=\varnothing$, $U_0=U$, $V_1 \cup V_2 \nsubseteq U$, and the start $y$ of $P$ lies in $U$, so the desired $w_1,v_1$ exist; also $J_0=J$, so $I_0$ exists by the definition of $J$. 
For $s>1$, if $V_1 \cup V_2 \subseteq U_{s-1}$, then by the process
\[
   V_3 \cup \cdots \cup V_{i+1} \subseteq I_{s-1} \setminus \{w_{s-1}\} \subseteq U_{s-1},
\]
and since $W_{s-1} \subseteq V_1 \cup V_2$,
\[
   U_{s-1} 
   = \bigl(V_1 \cup V_2 \setminus W_{s-1}\bigr) \cup W'_{s-1} \cup V_3 \cup \cdots \cup V_{i+1}
   = \bigl(J \setminus W_{s-1}\bigr) \cup W'_{s-1}
   = J_{s-1},
\]
so the algorithm should have terminated earlier, a contradiction. 
Hence $V_1 \cup V_2 \nsubseteq U_{s-1}\cup W_{s-1}$. 
Since the start $y$ of $P$ lies in $U_{s-1}\cup W_{s-1}$ (either still in $U_{s-1}$ or already in $W_{s-1}$), there exists an edge $\overrightarrow{w_s v_s}$ with $w_s \in U_{s-1}$ and $v_s \notin U_{s-1}$. 
As $\overrightarrow{w_s v_s}$ is an edge of $P$, we have $w_s,v_s \in J$. 
By the definition of $J$, there is an $(i+1)$-set $I_{s-1} \subseteq J$ containing $w_s,v_s$ with $F_{I_{s-1}}$ $(i,j,\varepsilon,X)$-balanced. 
Since points in $V_3,\dots,V_{i+1}$ are never modified, $V_3 \cup \cdots \cup V_{i+1} \subseteq J_{s-1}$, hence $I_{s-1} \subseteq J_{s-1} \cap J$. 
Therefore the process is well defined.

When the process ends with $U_s=J_s$, we set $g_u(J)=J_s$. 
By the construction,
\[
   F_{J_s}=F_{J_{s-1}}=\cdots=F_{J_0}=F_J
   \quad \text{and} \quad
   F_{\,J_s \setminus \{y_s\}}=\cdots =F_{\,J_0 \setminus \{y_0\}}=F_{\,J \setminus \{y\}}.
\]
Thus $g_u$ satisfies \ref{A1} and \ref{A2}.

It remains to prove \ref{A3}. 
Fix $\widetilde{J} \in g_u(\mathcal{J}_{u,i,j})$. 
If $u$ is supported by some preimage, then we do not change that preimage, so $\widetilde{J}$ itself can occur.

If Case~1 occurs, there are at most $(d+1)$ choices for $w'$, at most $\binom{d+1}{i+1}$ choices for $F_I$ (still determined by some $(i+1)$-set of $J'$), and in the bi-$(3d^2,6d^2)$ graph $G_{F_I}$ there are at most $6d^2$ choices of $w \in A_{F_I}$ adjacent to $w'$. 
Hence this case contributes at most
\[
   (d+1)\binom{d+1}{i+1}\cdot 6d^2 \;\le\; (6d)^{3d}.
\]

If Case~2 occurs, since the final $W'_s \subseteq J_s$, there are at most $2^{j+1}\le 2^{d+1}$ choices for $W'_s$, and at most $|W'_s|!\le (d+1)!$ possible orders of modification. 
Tracing the last modification, there exists an $(i+1)$-set $I'$ of $U_{|W'_s|-1}$ containing the last-added vertex such that the corresponding edge lies in $G_{F_{I'}}$, giving at most $\binom{j+1}{i+1}\cdot 6d^2 \le 6d^{\,d+2}$ choices at that step; similar bounds hold for each element of $W_s$. 
Therefore this case contributes at most
\[
   2^{d+1}\cdot (d+1)!\cdot 6d^{\,d+2} \;\le\; (12d)^{4d}
\]
preimages.

In total, there are at most 
\[
   1 + (6d)^{3d} + (12d)^{4d} \;\le\; (12d)^{5d}
\]
possible preimages of $\widetilde{J}$, proving \ref{A3}.\qedhere
\end{proof}

\subsection{Balanced supersaturation}

Now we are ready to show balanced supersaturation. Recall that $\mathcal{J}_{i,j}$ is the set of $(i,\varepsilon, X)$-good $(j+1)$-sets and we say an IGP $(j+1)$-set $J$ is \textit{$(i,\varepsilon,X)$-good} if $|X\cap F_J|\le nq^{-\frac{j}{d}}$ and there exists a partition of $(j+1)$-set with \[
  |V_1|=\left\lceil\tfrac{j-i}{2}\right\rceil+1,\qquad
  |V_2|=\left\lfloor\tfrac{j-i}{2}\right\rfloor+1,\qquad
  |V_3|=\dots=|V_{i+1}|=1,
\] such that for any $(i+1)$-set $I$ with $|I\cap V_s|=1$ for $s\in [i+1]$, the flat $F_I$ is $(i,j,\varepsilon ,X)$-balanced.

\begin{lemma}\label{low}
    There exists a constant $C$ such that the following statement holds. Let $X$ be a point set of $\mathbb{F}_q^d$ with $|X|=n\ge Cq$. Suppose that for some $1\le i<j\le d-1$ and $1/q\ll \varepsilon\ll 1/d$, there are at least $\Omega(n^{i+1})$ IGP $(i+1)$-sets $I$ in $X^{i+1}$ with $F_I$  $(i,j,\varepsilon,X)$-balanced and there are $\Omega(n^{j+1})$ $(i,\varepsilon,X)$-good IGP $(j+1)$-sets. Then there exists a collection $\mathcal{S}$ of coplanar sets in $X$ such that 
     \begin{enumerate}[label=\rm{(\roman*)}]
        \item\label{equation3} $|\mathcal{S}|=\Omega\l(n^{d+1}/q\r)$ and
        \item\label{equation4} $\Delta_s(\mathcal{S})\le O\l(n^{d+1-s}\cdot{q^{-\frac{d+1-s}{d}}}\r)$ for $s\in [d+1]$. 
    \end{enumerate}
\end{lemma}
\begin{proof}
    Let $\mathcal{I}$ be the family of IGP $(i+1)$-sets in $X$ with $F_I$ $(i,j,\varepsilon, X)$-balanced. Recall that $\mathcal{J}_{i,j}$ is the collection of $(i,\varepsilon, X)$-good $(j+1)$-sets. By the definition of $(i,j, \varepsilon, X)$-balanced, each such $F_I$ satisfies that $|F_I\cap X|\ge \varepsilon^{d-i}nq^{-\frac{j+1}{d}}.$ Note that for each $J$, $|V_1|=\l\lceil\frac{j-i}{2}\r\rceil+1\ge 2$. We define a root function $f:\mathcal{J}_{i,j}\rightarrow X$ so that for each $J\in \mathcal{J}_{i,j}$, $f(J)$ lies in the corresponding $V_1$ of~$J$.
 
We proceed in the following steps to construct $\mathcal{S}$.

\medskip
\textbf{Step 1.}  Choose arbitrarily an $(i,\varepsilon,X)$-good IGP $(j+1)$-set $J\in \mathcal{J}_{i,j}$. Write $J=\{y,x_1,\dots x_{j}\}$ with $y=f(J)$.

\textbf{Step 2.} Choose arbitrarily $u\in (F_J\setminus F_{J\setminus \{y\}})\cap X$. Then  $J\in \mathcal{J}_{u,i,j}$ by definition.

\textbf{Step 3.} Let $g_u:\mathcal{J}_{u,i,j}\rightarrow \binom{X}{j+1}$ be the map given by Lemma~\ref{map lemma}. Set $U:=g_u(J)\cup \{u\}$. By Lemma~\ref{map lemma}~\ref{A1}, $U$ is a critical coplanar $(j+2)$-set.

\textbf{Step 4.} Choose arbitrarily the remaining $(d-j-1)$ points $W=\{w_1,\dots ,w_{d-j-1}\}\subseteq X\setminus U$.

\textbf{Step 5.} Keep the coplanar $U\cup W$ with probability 
\[
  p \;=\; \frac{\varepsilon^{\,d} n}{\,q\,\bigl|(F_J\setminus F_{\,J\setminus\{y\}})\cap X\bigr|}\,.
\]
\medskip

    To see that these steps are feasible, it suffices to show that $p\le 1$. Let $I$ be an $(i,j,\varepsilon,X)$-balanced set of $J$ containing $y$. Since $y\notin F_{I\setminus \{y\}}$, we have $F_{I}\cap F_{J\setminus \{y\}}\neq F_{I}$, so $F_{I}\cap F_{J\setminus \{y\}}$ is contained in some $\l(i-1\r)$-flat of $F_{I}$. By the $(i,j,\varepsilon,X)$-balanced property, 
    \[|\l(F_I\setminus F_{J\setminus \{y\}}\r)\cap X|\;=\;|F_I\cap X|-|\l(F_I\cap F_{J\setminus \{y\}}\r)\cap X|\;\ge\; (1-\varepsilon)|F_I\cap X|\;\ge\; \varepsilon^{d}\,n\,q^{-\frac{j+1}{d}}.\]
    Since $F_{I}\subseteq F_J$, we obtain
   \[
  p \;=\; \frac{\varepsilon^{\,d} n}{\,q\,\bigl|(F_J\setminus F_{\,J\setminus\{y\}})\cap X\bigr|}
  \;\le\; \frac{\varepsilon^{\,d} n}{\,q\,\bigl|(F_I\setminus F_{\,J\setminus\{y\}})\cap X\bigr|}
  \;<\; 1.
\]
    By Lemma~\ref{blanceprop}, there exists $\Omega(n^{j+1})$ elements in $\mathcal{J}_{i,j}$. By Lemma~\ref{map lemma}~\ref{A3}, for each $u\in X$ and $J\in \mathcal{J}_{i,j,u}$, we have $|g_{u}^{-1}(J)|\le (12d)^{5d}$. Hence 

    \begin{align*}
        \mathbb{E}[|\mathcal S|]&\ge |\mathcal{J}_{i,j}|\cdot|\l(F_{J}\setminus F_{J\setminus\{y\}}\r)\cap X|\cdot (12d)^{-5d} \cdot (n-d)^{d-j-1}\cdot\frac{\varepsilon^dn}{q|\l(F_{J}\setminus F_{J\setminus\{y\}}\r)\cap X|}\frac{1}{d!}\\
        &\ge |\mathcal{J}_{i,j}|\cdot (n/2)^{d-j-1}\cdot (12d)^{-5d} \cdot \frac{\varepsilon^d n}{q}\\
        &\ge \Omega\l(\frac{n^{d+1}}{q}\r).
    \end{align*}

Since $|\mathcal{S}|$ is a sum of independent Bernoulli random variables, the Chernoff bound (\Cref{lem: Chernoff}) gives
$$\Pr\left[\,|\mathcal S|\le \mathbb{E}[|\mathcal S|]/2\,\right] \;\le\; \exp(-\Omega(\mathbb{E}[\mathcal{S}]))\;\le\; o(1).$$
Hence $|\mathcal{S}|\ge \Omega\l(\frac{n^{d+1}}{q}\r)$ with probability $1-o(1)$.

 Now we check~\ref{equation4}. We first bound expectations. To this end, we will only show that, for any $d$-set $D$ in $X$,
    \begin{align}\label{Delta34}
         \E[\deg(D)]\;\le\; O\l(n\,q^{-\frac{d-j-1}{d}}\r),
    \end{align}
   for any $(j+1)$-set $D$ in $X$,
    \begin{align}\label{Delta35}
        \E[\deg(D)]\;\le\; O\l(n^{d-j-1}q^{-\frac{d-1}{d}}\r),
    \end{align}
    for any $1$-set $D$ in $X$, 
    \begin{align}\label{Delta36}
        \E[\deg(D)]\;\le\; O(nq^{-1})
    \end{align}

     Note that, for any $1\le s\le s'\le d$, any $s$-set $D$ satisfies \[
  \deg(D)\;\le\; n^{\,s'-s}\,\max\{\,\deg(D') \;:\; D\subset D',\ |D'|=s'\,\}.
\] Thus~(\ref{Delta34}),~(\ref{Delta35}) and (\ref{Delta36}) would imply, for any $1\le s\le d$ and any $s$-set $D$,
     \begin{align}\label{DeltaExpectation1}
     \E[\deg(D)]\;\le\; O(n^{d+1-s}q^{-\frac{d+1-s}{d}}).
     \end{align}
     
     First, to see~(\ref{Delta34}), note that for a fixed $d$-set $D$, the remaining point have at most $n$ choices, and in step $4$, we keep the $(d+1)$-set with probability $\frac{\varepsilon^{d-i} n}{q|F_J\setminus F_{J\setminus \{y\}}\cap X|}$. Recall that we have 
     \[|F_J\setminus F_{J\setminus \{y\}}\cap X|\;\ge\; \varepsilon^{d}\,n\,q^{-\frac{j+1}{d}}.\] Thus the expectation of $\deg(D)$ can be estimated as follows.
\[
   \mathbb{E}[\deg(D)]
   \;\le\;
   n \cdot \frac{\varepsilon^{\,d}\, n}
      {\,q\,\bigl|(F_J \setminus F_{\,J \setminus \{y\}})\cap X\bigr|}
   \;\le\;
   \frac{n}{q^{(d-j-1)/d}}.
\]
    To see~(\ref{Delta35}), we fixed a $(j+1)$-set $D$. We need to count the number of sets $U\cup W$ generated by the steps above that contains $D$. First we decide which $j+1$ points in $U\cup W$ is specified by $D$ --- there are $\binom{d+1}{j+1}$ ways to do this. Note that the critical coplanar $(j+2)$ set $U=g_u(J)\cup \{u\}$ cannot be contained in $D$ (since $j+2>j+1$). Thus there exists a point $v\in U$ not specified. Next we specify all other $d-j-1$ unspecified points except $v$ --- there are at most $n^{d-j-1}$ choices. Since $U$ is critical coplanar, after all $i+1$ points except $v$ in $U$ have been specified, we have specified the $i$-flat $F$ which must contain $v$. By Lemma~\ref{map lemma}~\ref{A2}, we see that $F_{g_u(J)}=F_J$, where $J$ is the set chosen in step $2$ corresponding to $g_u(J)\cup \{u\}$. So there are at most $|F_{J}\cap X|$ ways to specify $v$.
    Note that $g_u(J)\cup u$ can be chosen at most $(12d)^{5d}$ times by Lemma~\ref{map lemma}~\ref{A3}.

    Since $y=f(J)$ is in the $V_1$ of $J$, by Lemma~\ref{key ob}, we always have 
    \[
   \frac{|F_J \cap X|}
        {\;\bigl|(F_J \setminus F_{\,J \setminus \{y\}})\cap X\bigr|}
   \;\le\; \varepsilon^{-d}\, q^{1/d}.
\]
    and so we can estimate the expectation of $\deg(D)$ as follows.
    \[
\mathbb{E}[\deg(D)]
\;\le\;
\binom{d+1}{j+1}\, n^{\,d-j-1}\,(12d)^{5d}\,|F_J\cap X|\,
\frac{\varepsilon^{\,d} n}{\,q\,\bigl|(F_J\setminus F_{\,J\setminus\{y\}})\cap X\bigr|}
\;\le\;
O\!\left(\frac{n^{\,d-j}}{q^{(d-1)/d}}\right).
\]
    To see~(\ref{Delta36}), we fixed a point $D=\{v\}$. Similar to the above, we count the number of sets $U\cup W$ generated by the steps that contains $D$. First we decide which point in $U\cup W$ is specified by $D$ --- there are $d+1$ ways to do this. 
    
    If $v$ plays the role of some vertex in $g_u(J)\cup \{w_1,\dots, w_{d-j-1}\}$, then there are at most $n^{d-1}$ ways to specify the other points in  $g_u(J)\cup \{w_1,\dots, w_{d-j-1}\}$. Fix a choice of $g_u(J)\cup \{w_1,\dots, w_{d-j-1}\}$. By Lemma~\ref{map lemma}, we can determine the $F_J$ and $F_{J\setminus f(J)}$.
    
    If $v$ plays the role of $u$, by Lemma~\ref{map lemma}~\ref{A2}, there exists a vertex $y'$ such that $F_{J\setminus \{f(J)\}}=F_{g_u(J)\setminus \{y'\}}$. Note that there are $j+1$ ways to specified the position of $y'$ in $g_u(J)$. For the other $d-1$ points except the $y'$, there are at most $n^{d-1}$ ways to specify them. Fix a choice of $g_{u}(J)\setminus\{y'\}\cup \{w_1,\dots, w_{d-j}\}\cup \{u\}$. Since $u$ is supported by $g_{u}(J)$, we see that $F_{g_u(J)\setminus \{y'\}\cup \{u\}}=F_J$. Moreover, the flat $F_{J\setminus \{f(J)\}}$ is $F_{g_{u}(J)\setminus \{y'\}}$.    
    
 Note that in both cases, we have determined the flat $F_J$ and $F_{J\setminus \{f(J)\}}$ by~Lemma~\ref{map lemma}~\ref{A2}. So there are $|F_J\setminus F_{J\setminus \{f(J)\}}|$ choices for the last vertex $u$ ($f(J)$, resp). Thus we can estimate the expectation of $\deg(D)$ as follows.
    $$
    \E[\deg(D)]\;\le\; (d+1)\cdot (j+1)\cdot n^{d-1}|(F_J\setminus F_{J\setminus\{y\}})\cap X|\cdot\frac{\varepsilon^{d} n}{q|(F_J\setminus F_{J\setminus \{f(J)\}})\cap X)|}\;\le\; O\l(\frac{n^{d-j}}{q}\r).
    $$

    Now we have inequalities~(\ref{DeltaExpectation1}). Thus for any $1\le s\le d$, there exists a sufficiently large constant $C_s$ such that for any $s$-set $D$ in $X$, we have
    $\E[\deg(D)]\le C_sn^{d+1-s}q^{-\frac{d+1-s}{d}}.$ For any fixed $s$-set $D$, note that $\deg(D)$ is a sum of independent Bernoulli random variables. Thus by the Chernoff bound (Lemma~\ref{lem: Chernoff}), together with the fact that $q\le O(n)$, $2\le d$ and $s\le d$, we have
    $$
    \Pr\l[\deg(D)\;\ge\; 20C_s\,n^{d+1-s}\,q^{-\frac{d+1-s}{d}}\r]\;\le\; 2\exp\l(-10C_s\,n^{d+1-s}\,q^{-\frac{d+1-s}{d}}\r)\;\le\; \exp\l(-\Omega\l(n^{\frac{d-1}{d}}\r)\r).
    $$

    By the union bound, for any $1\le s\le d$,
    $$
    \Pr\l[\Delta_{s}(\mathcal{S})\;\ge\; 20C_s\,n^{d+1-s}\,q^{-\frac{d+1-s}{d}}\r]\;\le\; \binom{n}{s}\exp\l(-\Omega\l(n^{\frac{d-1}{d}}\r)\r)\;\le\; o(1).
    $$
     
     Hence, with probability $1-o(1)$, we have $\Delta_s(\mathcal{S})\le O\l(n^{d+1-s}q^{-\frac{d+1-s}{d}}\r)$.

    Therefore, we can find the collection $\mathcal{S}$ as desired.\qedhere
\end{proof}
\section{Proof of Theorem~\ref{main}}\label{proof1}
\begin{proof}[Proof of Theorem~\ref{main}]
    By Lemma~\ref{concentrationprop}, one of the following holds.

    \textbf{Case 1.} There are $\Omega(n^{d+1})$ coplanar sets.  
    Then by Lemma~\ref{struc}, there exists a collection~$\mathcal{S}$ as desired.

    \textbf{Case 2.} There are $\Omega(n^{i+1})$ $(i,j,\varepsilon,X)$-balanced sets with $j\le i$ for some $\varepsilon>0$.  
    Then by Lemma~\ref{high}, there exists a collection~$\mathcal{S}$ as desired.

    \textbf{Case 3.} There are $\Omega(n^{i+1})$ $(i,j,\varepsilon,X)$-balanced sets with $j> i$ for some $\varepsilon>0$.  
    Then Lemma~\ref{concentrationprop} gives $\Omega(n^{j+1})$ $(i,\varepsilon,X)$-good $(j+1)$-sets, and by Lemma~\ref{low}, there exists a collection~$\mathcal{S}$ as desired. 

    In all cases, the conclusion follows. 
\end{proof}

\section{Concluding remarks}\label{concluding}
\begin{itemize}
    \item Let $N(q,d,k)$ denote the number of IGP $k$-sets in $\F^d_q$. In the case $d=2$, Roche-Newton and Warren~\cite{roche2022arcs} discovered a threshold of $k$ where the behavior of $N(q,2,k)$ suddenly change:  when $k\ll q^{1/2}$, $N(q,2,k)=\binom{\Theta(q^2)}{k}$, whereas when $k\gg q^{1/2+o(1)}$, $N(q,2,k)\le\binom{O(q^{3/2+o(1)})}{k}$. The second range was later improved by the first author, Liu, the second author, and Zeng~\cite{CLNZ}, who showed that when $k\gg q^{1/2+o(1)}$, $N(q,2,k)\le\binom{\Theta(q)}{k}$. More recently, Nenadov~\cite{nenadov2024number} made a further improvement and proved that when $k\gg q^{1/2+o(1)}$, $N(q,2,k)\le\binom{(1+o(1))q}{k}$. Using \Cref{main} we can show that following result for $d\ge 3$.
    \begin{theorem}
    For fixed integer $d\ge 3$ and all sufficiently large prime powers $q$,  we have
    $$
    N(q,d,k)\;=\;\binom{\Theta(q^d)}{k}\quad \text{when } k\ll q^{1/d},
    $$
    and 
    $$
    N(q,d,k)\;=\;\binom{\Theta(q)}{k}\quad\text{when }k\ge q^{1/d+o(1)}.
    $$
    \end{theorem}

    The proof is a straightforward generalization of the proofs of Theorem 1 in~\cite{roche2022arcs} and Theorem 1.5 in~\cite{CLNZ}, so we omit the details. Inspired by the work of Nenadov for $d=2$, we propose the following conjecture.
    \begin{conjecture}
    For any fixed integer $d\ge 3$ and sufficiently large prime powers $q$, when $k\ge q^{1/d+o(1)}$, we have
    $$
    N(q,d,k)\;=\;\binom{(1+o(1))q}{k}.
    $$    
    \end{conjecture}
    
    \item A point set $X$ in $\mathbb{F}_q^d$ is said to be \textit{$(k,c)$-evasive} if the intersection between $X$ and any $k$-flat has cardinality less than $c$. In this language, an IGP set is just a $(d-1,d+1)$-evasive set. Evasive sets arise in various areas of mathematics, including Ramsey theory~\cite{pudlak2004pseudorandom}, incidence geometry~\cite{sudakov2022evasive,milojevic2024incidence}, error-correcting codes~\cite{guruswami2011linear,dvir2012subspace}, and extremal graph theory~\cite{blagojevic2013turan,bukh2024extremal}.

    Recently, Lim, the second author, and Zeng~\cite{lim2025evasive} proved that the number of $(k,c)$-evasive sets is $2^{\Theta(q^{\,d-k})}$, where $c$ is a fixed integer sufficiently large relative to $d$ and $k$, and $q$ is a sufficiently large prime power. In their proof they use a balanced supersaturation result for evasive sets (Lemma 5.1 in~\cite{lim2025evasive}), which is suboptimal but sufficient for their purpose. We believe the techniques  developed in the present paper could be generalized to yield, in an appropriate sense, ``optimal'' balanced supersaturation for evasive sets, which in turn could lead to ``optimal'' fixed-size counting results and random Turán results for evasive sets.

\end{itemize}

\bibliography{ref}
\end{document}